\documentclass[11pt]{article}
\usepackage{amsthm,amssymb,amsmath,amscd}
\usepackage[all]{xy}
\usepackage{mathrsfs}
\usepackage{multirow}

\usepackage{mathptm}

\newtheorem{Def}{Definition}[section]
\newtheorem{Prop}[Def]{Proposition}
\newtheorem{Theo}[Def]{Theorem}
\newtheorem{Lem}[Def]{Lemma}
\newtheorem{Koro}[Def]{Corollary}

\newcommand{\E}[4]{{\rm E}_{#1}^{#2}(#3, #4)}
    \newcommand{\Ex}[3]{{\rm E}_{#1}^{#2}(#3)}
    \newcommand{\EP}[4]{{\rm E}_{#1}^{#2\bullet}(#3, #4)}

\newcommand{\add}{{\rm add}}
\newcommand{\con}{{\rm con}}

\newcommand{\Hom}{{\rm Hom }}
\newcommand{\rad}{{\rm rad}}
\newcommand{\soc}{{\rm soc}}
\renewcommand{\top}{{\rm top}}

\newcommand{\End}{{\rm End}}
\newcommand{\Ext}{{\rm Ext}}
\newcommand{\Coker}{{\rm Coker}}
\newcommand{\Ker}{{\rm Ker}}
\newcommand{\cpx}[1]{#1^{\bullet}}
\newcommand{\D}[1]{\mathscr{D}(#1)}

\newcommand{\Db}[1]{ \mathscr{D}^{\rm b}(#1)}
\newcommand{\C}[1]{\mathscr{C}(#1)}

\newcommand{\Cb}[1]{{\mathscr{C}^b}(#1)}
\newcommand{\K}[1]{\mathscr{K}(#1)}

\newcommand{\Kf}[1]{\mathscr{K}^-(#1)}
\newcommand{\Kb}[1]{ \mathscr{K}^{\rm b}(#1)}
\newcommand{\modcat}[1]{#1\mbox{{\rm -mod}}}
\newcommand{\stmodcat}[1]{#1\mbox{{\rm -{\underline{mod}}}}}
\newcommand{\pmodcat}[1]{#1\mbox{{\rm -proj}}}
\newcommand{\imodcat}[1]{#1\mbox{{\rm -inj}}}
\newcommand{\opp}{^{\rm op}}
\newcommand{\otimesL}{\otimes^{\rm\bf L}}

\newcommand{\HomP}{{\rm Hom}^{\bullet}}

\newcommand{\lra}{\longrightarrow}
\newcommand{\ra}{\rightarrow}
\newcommand{\lraf}[1]{\stackrel{#1}{\lra}}

\begin{document}

{\Large \bf
\begin{center}
Derived equivalences for $\Phi$-Auslander-Yoneda algebras

\end{center}}
\medskip

\centerline{{\bf Wei Hu}, {\bf Changchang Xi$^*$}}
\begin{center} School of Mathematical Sciences; \\
Laboratory of Mathematics and Complex Systems, \\ Beijing Normal University, 100875 Beijing,\\
 People's Republic of  China \\E-mail: huwei@bnu.edu.cn \quad  xicc@bnu.edu.cn\\
\end{center}

\begin{abstract}
In this paper, we introduce $\Phi$-Auslander-Yoneda algebras in a
triangulated category with $\Phi$ a parameter set in $\mathbb N$,
and present a method to construct new derived equivalences between
these $\Phi$-Auslander-Yoneda algebras (not necessarily Artin
algebras), or their quotient algebras, from a given almost
$\nu$-stable derived equivalence. As consequences of our method, we
have: (1) Suppose that $A$ and $B$ are representation-finite,
self-injective Artin algebras with $_AX$ and $_BY$ additive
generators for $A$ and $B$, respectively. If $A$ and $B$ are
derived-equivalent, then the $\Phi$-Auslander-Yoneda algebras of $X$
and $Y$ are derived-equivalent for every admissible set $\Phi$. In
particular, the Auslander algebras of $A$ and $B$ are both
derived-equivalent and stably equivalent. (2) For a self-injective
Artin algebra $A$ and an $A$-module $X$, the $\Phi$-Auslander-Yoneda
algebras of $A\oplus X$ and $A\oplus \Omega_A(X)$ are
derived-equivalent for every admissible set $\Phi$, where $\Omega$
is the Heller loop operator. Motivated by these derived equivalences
between $\Phi$-Auslander-Yoneda algebras, we consider constructions
of derived equivalences for quotient algebras, and show, among other
results, that a derived equivalence between two basic self-injective
algebras may transfer to a derived equivalence between their
quotient algebras obtained by factorizing out socles.
\end{abstract}

\renewcommand{\thefootnote}{\alph{footnote}}
\setcounter{footnote}{-1} \footnote{ $^*$ Corresponding author.
Email: xicc@bnu.edu.cn; Fax: 0086 10 58802136; Tel.: 0086 10
58808877.}
\renewcommand{\thefootnote}{\alph{footnote}}
\setcounter{footnote}{-1} 
\renewcommand{\thefootnote}{\alph{footnote}}
\setcounter{footnote}{-1} 

\section{Introduction}
Derived categories and derived equivalences were introduced by
Grothendieck and Verdier in \cite{v}. As is known, they have widely
been used in many branches of mathematics and physics. One of the
fundamental problems in the study of derived categories and derived
equivalences is: how to construct derived equivalences ? On the one
hand, Rickard's beautiful Morita theory for derived categories can
be used to find all rings that are derived-equivalent to a given
ring $A$ by determining all tilting complexes over $A$ (see
\cite{RickMoritaTh} and \cite{RickDFun}). On the other hand, a
natural course of investigation on derived equivalences should be
constructing new derived equivalences from given ones. In this
direction, Rickard used tensor products and trivial extensions to
produce new derived-equivalences in \cite{RickMoritaTh, Rickard3},
Barot and Lenzing employed one-point extensions to transfer certain
a derived equivalence to a new one in \cite{BL}. Up to now, however,
it seems that not much is known for constructing new derived
equivalences based on given ones.

In this paper, we continue the consideration in this direction and
provide, roughly speaking, two methods to construct new derived
equivalences from given ones. One is to form $\Phi$-Auslander-Yoneda
algebras (see Section \ref{subsect3.1} for definition) of
generators, or cogenerators over derived-equivalent algebras, and
the other is to form quotient algebras of derived-equivalent
algebras. We point out that our family of $\Phi$-Auslander-Yoneda
algebras include Auslander algebras, generalized Yoneda algebras and
some of their quotients. Thus our method produces also derived
equivalences between infinite-dimensional algebras.

To state our results, we first introduce a few terminologies.

Suppose that $F$ is a derived equivalence between two Artin algebras
$A$ and $B$, with the quasi-inverse functor $G$. Further, suppose
that
$$\cpx{T}: \quad \cdots \lra 0\lra T^{-n}\lra
\cdots \lra T^{-1}\lra T^0\lra 0\lra \cdots$$ is a radical tilting
complex over $A$ associated to $F$, and suppose that
$$\cpx{\bar{T}}: \quad \cdots \lra 0\lra \bar{T}^{0}
\lra \bar{T}^{1}\lra \cdots \lra \bar{T}^n\lra 0\lra \cdots
$$ is a radical tilting complex over $B$ associated to $G$. The functor $F$ is called {\it almost
$\nu$-stable} if $\add(\bigoplus_{i=-1}^{-n} T^i)$ =
$\add(\bigoplus_{i=-1}^{-n} \nu_A T^i)$, and
$\add(\bigoplus_{i=1}^{n} \bar{T}^i)$ = $\add(\bigoplus_{i=1}^n
\nu_B\bar{T}^i)$, where $\nu_A$ is the Nakayama functor for $A$. We
have shown in \cite{HX3} that an almost $\nu$-stable functor $F$
induces an equivalence functor $\bar{F}$ between the stable module
categories $\stmodcat{A}$ and $\stmodcat{B}$. For further
information on almost $\nu$-stable derived equivalences, we refer
the reader to \cite{HX3}.

For a module $M$ over an algebra $A$, the \emph{generalized Yoneda
algebra} of $M$ is defined by Ext$^*_A(M):=\bigoplus_{i\ge
0}\Ext^i_A(M,M)$. In case $M=A/\rad(A)$, the algebra $\Ext^*_A(M)$
is called the \emph{Yoneda algebra} of $A$ in literature. We shall
extend this notion to a more general situation, and introduce the
$\Phi$-Auslander-Yoneda algebras with $\Phi$ a parameter set in
$\mathbb N$ (for details see Subsection \ref{subsect3.1} below). We
notice that a $\Phi$-Auslander-Yoneda algebra may not be an Artin
algebra in general.

Our main result on $\Phi$-Auslander-Yoneda algebras of modules reads
as follows:

\begin{Theo}
Let $A$ and $B$ be two Artin algebras, and let $\bar{F}:
\stmodcat{A}\lra\stmodcat{B}$ be the stable equivalence induced by
an almost $\nu$-stable derived equivalence $F$ between $A$ and $B$.
Suppose that $X$ is an $A$-module, we set $M:=A\oplus X$ and
$N:=B\oplus \bar{F}(X)$. Let $\Phi$ be an admissible subset of
$\mathbb{N}$, and define the $\Phi$-Auslander-Yoneda algebra of $M$
as $\Ex{A}{\Phi}{M}:=\bigoplus_{i\in \Phi}\emph{Ext}^i_A(M,M)$.
Then:

\smallskip
$(1)$ The $\Phi$-Auslander-Yoneda algebras $\Ex{A}{\Phi}{M}$ and
$\Ex{B}{\Phi}{N}$ are derived-equivalent.

$(2)$ If $\Phi$ is finite, then there is an almost $\nu$-stable
derived equivalence between $\Ex{A}{\Phi}{M}$ and $\Ex{B}{\Phi}{N}$.
Thus $\Ex{A}{\Phi}{M}$ and $\Ex{B}{\Phi}{N}$ are also stably
equivalent. In particular, there is an almost $\nu$-stable derived
equivalence and a stable equivalence between $\End_A(M)$ and
$\End_B(N)$.\label{theo1}
\end{Theo}

A dual version of Theorem \ref{theo1} can be seen in Corollary
\ref{dual} below.

Since Auslander algebra and generalized Yoneda algebra are two
special cases of $\Phi$-Auslander-Yoneda algebras, Theorem
\ref{theo1} provides a large variety of derived equivalences between
Auslander algebras, and between generalized Yoneda algebras, or
their quotient algebras. Note that Theorem \ref{theo1} (2) extends a
result in \cite[Proposition 6.1]{HX3}, where algebras were assumed
to be finite-dimensional over a field, in order to employ two-sided
tilting complexes in proofs, and where only endomorphism algebras
were considered instead of general Auslander-Yoneda algebras. The
existence of two-sided tilting complexes is guaranteed for Artin
$R$-algebras that are projective as $R$-modules \cite{RickDFun}. For
general Artin algebras, however, we do not know the existence of
two-sided tilting complexes. Hence, in this paper, we have to
provide a completely different proof to the general result, Theorem
\ref{theo1}.

As a direct consequence of Theorem \ref{theo1}, we have the
following corollary concerning the Auslander algebras of
self-injective algebras.

\begin{Koro}
$(1)$  For a self-injective Artin algebra $A$ and an $A$-module $Y$,
the $\Phi$-Auslander-Yoneda algebras $\Ex{A}{\Phi}{A\oplus Y}$ and
$\Ex{A}{\Phi}{A\oplus \Omega_A(Y)}$ are derived-equivalent, where
$\Omega$ is the Heller loop operator.

$(2)$ Suppose that $A$ and $B$ are self-injective Artin algebras of
finite representation type with $_AX$ and $_BY$ additive generators
for {\rm $A$-mod} and {\rm $B$-mod}, respectively. If $A$ and $B$
are derived-equivalent, then

$(i)$ the Auslander algebras of $A$ and $B$ are both derived and
stably equivalent.

$(ii)$ The generalized Yoneda algebras $\emph{Ext}^*_A(X)$ and
$\emph{Ext}^*_B(Y)$ of $X$ and $Y$ are
derived-equivalent.\label{1.2}
\end{Koro}

Notice that, in Corollary \ref{1.2}, the Auslander algebra of $A$ is
a quotient algebra of the generalized Yoneda algebra $\Ext_A^*(X)$
of the additive generator $X$. The next result shows another way to
construct derived equivalences for quotient algebras.

\begin{Theo}   Let $F:\Db{A}\longrightarrow\Db{B}$ be a derived equivalence between
two self-injective basic Artin algebras $A$ and $B$. Suppose that
$P$ is a direct summand of $_AA$, and $Q$ is a direct summand of
$_BB$ such that $F(\soc(P))$ is isomorphic to $\soc(Q)$, where
$\soc(P)$ denotes the socle of the module $P$. Then the quotient
algebras $A/\soc(P)$ and $B/\soc(Q)$ are derived-equivalent.
\label{theo2}
\end{Theo}

The structure of this paper is organized as follows. In Section
\ref{sect2}, we make some preparations for later proofs.  In Section
\ref{sect3}, we introduce the $\Phi$-Auslander-Yoneda algebras and
prove Theorem \ref{theo1} and its dual version, Corollary
\ref{dual}, which produces derived equivalences between the
endomorphism algebras of cogenerators. Furthermore, we deduce a
series of consequences of Theorem \ref{theo1} for self-injective
algebras, including Corollary \ref{1.2}. In Section \ref{sect4}, we
provide several methods to construct derived equivalences for
quotient algebras. First, we give a general criterion, and then
apply it to self-injective algebras modulo socles, and to algebras
modulo annihilators. In particular, we show Theorem \ref{theo2}, and
point out a class of derived equivalences satisfying the conditions
in Theorem \ref{theo2}.

\section{Preliminaries \label{sect2}}

In this section, we shall recall basic definitions and facts on
derived categories and derived equivalences, which are elementary
elements in our proofs.

Throughout this paper, $R$ is a fixed commutative Artin ring. Given
an $R$-algebra $A$, by an $A$-module we mean a unitary left
$A$-module; the category of all finitely generated $A$-modules is
denoted by $\modcat{A}$, the full subcategory of $A$-mod consisting
of projective (respectively, injective) modules is denoted by
$\pmodcat{A}$ (respectively, $\imodcat{A}$). The stable module
category $\stmodcat{A}$ of $A$ is, by definition, the quotient
category of $\modcat{A}$ modulo the ideal generated by homomorphisms
factorizing through projective modules. An equivalence between the
stable module categories of two algebras is called a {\em stable
equivalence}

An $R$-algebra $A$ is called an \emph{Artin $R$-algebra} if $A$ is
finitely generated as an $R$-module. For an Artin $R$-algebra $A$,
we denote by $D$ the usual duality on $\modcat{A}$, and by $\nu_A$
the Nakayama functor $D\Hom_A(-, {}_AA): A\mbox{-proj}\lra
A\mbox{-inj}$. For an $A$-module $M$, we denote by $\Omega_A(M)$ the
first syzygy of $M$, and call $\Omega_A$ the \emph{Heller loop
operator} of $A$. In this paper, we mainly concentrate us on Artin
algebras and finitely generated modules.

Let $\cal C$ be an additive category.

For two morphisms $f:X\rightarrow Y$ and $g:Y\rightarrow Z$ in $\cal
C$, we write $fg$ for their composition. But for two functors
$F:\mathcal{C}\rightarrow \mathcal{D}$ and
$G:\mathcal{D}\rightarrow\mathcal{E}$ of categories, we write $GF$
for their composition instead of $FG$. For an object $X$ in
$\mathcal{C}$, we denote by $\add(X)$ the full subcategory of $\cal
C$ consisting of all direct summands of finite direct sums of copies
of $X$. An object $X$ in $\cal C$ is called an \emph{additive
generator} for $\cal C$ if add$(X)={\cal C}$.

By a complex $\cpx{X}$ over $\cal C$ we mean a sequence of morphisms
$d_X^{i}$ between objects $X^i$ in $\cal C$: $ \cdots \rightarrow
X^i\stackrel{d_X^i}{\lra}
X^{i+1}\stackrel{d_X^{i+1}}{\lra}X^{i+2}\rightarrow\cdots, $ such
that $d_X^id_X^{i+1}=0$ for all $i \in {\mathbb Z}$, and write
$\cpx{X}=(X^i, d_X^i)$. For a complex $\cpx{X}$, the {\em brutal
truncation} $\sigma_{<i}\cpx{X}$ of $\cpx{X}$ is a subcomplex of
$\cpx{X}$ such that $(\sigma_{<i}\cpx{X})^k$ is $X^k$ for all $k<i$
and zero otherwise. Similarly, we define $\sigma_{\geqslant
i}\cpx{X}$. For a fixed $n\in {\mathbb Z}$, we denote by
$\cpx{X}[n]$ the complex obtained from $\cpx{X}$ by shifting $n$
degrees, that is, $(\cpx{X}[n])^0=X^n$.

The category of all complexes over $\cal C$ with chain maps is
denoted by $\C{\mathcal C}$. The homotopy category of complexes over
$\mathcal{C}$ is denoted by $\K{\mathcal C}$. When $\cal C$ is an
abelian category, the derived category of complexes over $\cal C$ is
denoted by $\D{\cal C}$. The full subcategory of $\K{\cal C}$ and
$\D{\cal C}$ consisting of bounded complexes over $\mathcal{C}$ is
denoted by $\Kb{\mathcal C}$ and $\Db{\mathcal C}$, respectively. As
usual, for an algebra $A$, we simply write $\C{A}$ for
$\C{\modcat{A}}$, $\K{A}$ for $\K{\modcat{A}}$ and $\Kb{A}$ for
$\Kb{\modcat{A}}$. Similarly, we write $\D{A}$ and $\Db{A}$ for
$\D{\modcat{A}}$ and $\Db{\modcat{A}}$, respectively.

It is well-known that, for an $R$-algebra $A$, the categories
$\K{A}$ and $\D{A}$ are triangulated categories. For basic results
on triangulated categories, we refer the reader to the excellent
books \cite{HappelTriangle} and \cite{neeman}.

Let $A$ be an Artin algebra. Recall that a homomorphism $f: X\ra Y$
of $A$-modules is called a \emph{radical map} if, for any module $Z$
and homomorphisms $h: Z\ra X$ and $g: Y\ra Z$, the composition $hfg$
is not an isomorphism. A complex over $\modcat{A}$ is called a
\emph{radical} complex if all of its differential maps are radical
maps. Every complex over $\modcat{A}$ is isomorphic  to a radical
complex in the homotopy category $\K{A}$. If two radical complexes
$\cpx{X}$ and $\cpx{Y}$  are isomorphic in $\K{A}$, then $\cpx{X}$
and $\cpx{Y}$ are isomorphic in $\C{A}$.

Two $R$-Artin algebras $A$ and $B$ are said to be
\emph{derived-equivalent} if their derived categories $\Db{A}$ and
$\Db{B}$ are equivalent as triangulated categories. By a result of
Rickard (see Lemma \ref{rickard} below), two algebras $A$ and $B$
are derived-equivalent if and only if $B$ is isomorphic to the
endomorphism algebra $\End_{\Kb{A}}(\cpx{T})$ of a tilting complex
$\cpx{T}$ over $A$. Recall that a complex $\cpx{T}$ in
$\Kb{\pmodcat{A}}$ is called a \emph{tilting complex} over $A$ if it
satisfies

(1) $\Hom_{\Kb{\pmodcat{A}}}(\cpx{T},\cpx{T}[n])=0$ for all $n\ne
0$, and

(2) $\add(\cpx{T})$ generates $\Kb{\pmodcat{A}}$ as a triangulated
category.

\smallskip
It is known that, given a derived equivalence $F$ between
$A$ and $B$, there is a unique (up to isomorphism in $\Kb{A}$)
tilting complex $\cpx{T}$ over $A$ such that $F(\cpx{T})\simeq B$.
This complex $\cpx{T}$ is called a tilting complex \emph{associated}
to $F$. Recall that a complex $\cpx{X}$ of $A$-modules is called
\emph{self-orthogonal} if $\Hom_{\Db{A}}(\cpx{X},\cpx{X}[i])=0$ for
every $i\ne 0$.

The following lemma, proved in \cite[Lemma 2.2]{HX3}, will be  used
frequently in our proofs below.

\begin{Lem}\label{kdiso}
Let $A$ be an arbitrary ring with identity, and let $A\mbox{\rm
-Mod}$ be the category of all left (not necessarily finitely
generated) $A$-modules. Suppose that $\cpx{X}$ is a complex over
$A\mbox{\rm -Mod}$ bounded above, and that $\cpx{Y}$ is a complex
over $A\mbox{\rm -Mod}$ bounded below. Let $m$ be an integer. If one
of the following two conditions holds:

$(1)$ $X^i$ is projective for all $i>m$ and $Y^j=0$ for all $j<m$,

$(2)$ $Y^j$ is injective for all $j<m$  and $X^i=0$ for all $i>m$,

\noindent then  the localization functor $\theta:$ $ \K{A\mbox{\rm
-Mod}}$ $\lra \D{A\mbox{\rm -Mod}}$ induces an isomorphism
 $\theta_{\cpx{X},\cpx{Y}}:$ \newline $ \Hom_{\K{A\mbox{\rm
-Mod}}}(\cpx{X},\cpx{Y})$ $\lra \Hom_{\D{A\mbox{\rm
-Mod}}}(\cpx{X},\cpx{Y})$.
\end{Lem}

Thus, for the complexes $\cpx{X}$ and $\cpx{Y}$ given in Lemma
\ref{kdiso}, the computation of morphisms from $\cpx{X}$ to
$\cpx{Y}$ in $\D{A\mbox{\rm -Mod}}$ is reduced to that in
$\K{A\mbox{\rm -Mod}}$.

For later reference, we cite the following fundamental result on
derived equivalences by Rickard (see \cite[Theorem
6.4]{RickMoritaTh}) as a lemma.

\begin{Lem}$\cite{RickMoritaTh}$\label{rickard}
Let $\Lambda$ and $\Gamma$ be two rings. The following conditions
are equivalent:

$(a)$ $\Kf{\Lambda\emph{-Proj}}$ and $\Kf{\Gamma\emph{-Proj}}$ are
equivalent as triangulated categories;

$(b)$$\Db{\Lambda\emph{-Mod}} and \Db{\Gamma\emph{-Mod}}$ are
equivalent as triangulated categories;

$(c)$ $\Kb{\Lambda\emph{-Proj}}$ and $\Kb{\Gamma\emph{-Proj}}$ are
equivalent as triangulated categories;

$(d)$ $\Kb{\Lambda\emph{-proj}}$ and $\Kb{\Gamma\emph{-proj}}$ are
equivalent as triangulated categories;

$(e)$ $\Gamma$ is isomorphic to $\End(T)$, where $T$ is a tilting
complex in $\Kb{\Lambda\emph{-proj}}$.

Here $\Lambda\emph{-Proj}$ stands for the full subcategory of
$\Lambda\emph{-Mod}$ consisting of all projective $\Lambda$-modules.
\end{Lem}

Two rings $\Lambda$ and $\Gamma$ are called
\emph{derived-equivalent} if one of the above conditions (a)-(e) is
satisfied. For Artin algebras, the two definitions of a derived
equivalence coincide with each other.

\section{Derived equivalences for $\Phi$-Auslander-Yoneda algebras \label{sect3}}
As is known, Auslander algebra is a key to characterizing
representation-finite algebras, and Yoneda algebra plays a role in
the study of the graded module theory of Koszul algebras. In this
section, we shall first unify the two notions and introduce the
so-called $\Phi$-Auslander-Yoneda algebra of an object in a
triangulated category, where $\Phi$ is a parameter subset of
$\mathbb N$, and then construct new derived equivalences between
these $\Phi$-Auslander-Yoneda algebras from a given almost
$\nu$-stable derived equivalence. In particular, Theorem \ref{theo1}
will be proved, and  a series of its consequences will be deduced in
this section.

\subsection{Admissible sets and Auslander-Yoneda algebras\label{subsect3.1}}

First, we introduce some special subsets of the set
$\mathbb{N}:=\{0,1,2,\cdots,\}$ of the natural numbers, and then
define a class of algebras called Auslander-Yoneda algebras.

A subset $\Phi$ of $\mathbb{N}$ containing $0$ is called an {\em
admissible subset} of $\mathbb{N}$ if the following condition is
satisfied:

{\it  If $i, j$ and $k$ are in $\Phi$ such that $i+j+k\in\Phi$, then
$i+j\in\Phi$ if and only if $j+k\in\Phi$.}

\medskip
For instance, the sets $\{0, 3, 4\}$, $\{0, 1, 2,3, 4\}$ are
admissible subsets of $\mathbb{N}$. The following is a family of
admissible subsets of $\mathbb N$.

Let $n$ be a positive integer, and let $m$ be a positive integer or
positive infinity. We define
$$\Phi(n, m):=\{xn \mid x\in\mathbb{N}, 0\leq x< m+1\}.$$
Then $\Phi(n,m)$ is  an admissible subset in ${\mathbb N}$. Clearly,
we have $\Phi(1,\infty)=\mathbb{N}, \Phi(1, 0)=\{0\}$, and $\Phi(1,
m)=\{0, 1, 2,\cdots, m\}$.

Admissible subsets of $\mathbb N$ have the following simple
properties.

\begin{Prop}
  $(1)$ If $\Phi$ is an admissible subset of $\mathbb{N}$, then so is
  $m\Phi:=\{mx\mid x\in\Phi\}$ for every $m\in\mathbb{N}$.

  $(2)$ If $\Phi_1$ and $\Phi_2$ are admissible subsets of
  $\mathbb{N}$, then so is $\Phi_1\cap\Phi_2$. Moreover, the
  intersection of a family of admissible subsets of $\mathbb N$ is
  admissible.

  $(3)$ For a subset $\Phi\subseteq\mathbb{N}$ with $0\in \Phi$, the set
  $\Phi^m:=\{x^m\mid x\in\Phi\}$ is an admissible subset of
  $\mathbb{N}$ for every integer $m\geq 3$.
\end{Prop}

{\it Proof.}
 The statements (1) and (2) follow easily from the definition of an
 admissible subset. Now we consider (3). We pick an integer $m\ge 3$. Let $i^m,j^m,k^m$ and $l^m$ be in
 $\Phi^m$ such that $i^m+j^m+k^m=l^m$. If $i^m+j^m\in\Phi^m$, then
 $i^m+j^m=t^m$ for some $t\in\Phi$. By Fermat's Last Theorem, one of the
 integers $i$ and $j$ is zero. If $j=0$, then $j^m+k^m=k^m\in\Phi^m$.
 If $i=0$, then $j^m+k^m=l^m\in \Phi^m$. Similarly, we can show that if
 $j^m+k^m\in\Phi^m$, then $i^m+j^m\in \Phi$. Hence the set
 $\Phi^m$ is an admissible subset of $\mathbb{N}$.
$\square$

Note that $\Phi^2$ is not necessarily admissible in $\mathbb{N}$
even if $\Phi$ is an admissible subset of $\mathbb{N}$. For
instance, if $\Phi=\{0, 3, 4, 5, 12, 13\}$, then $\Phi$ is
admissible. Clearly, $3^2+4^2+12^2=13^2\in\Phi^2$ and
$3^2+4^2=5^2\in\Phi^2$, but $4^2+12^2\not\in\Phi^2$, so $\Phi^2$ is
not admissible.

\medskip
Now, we use admissible subsets to define a class of associative
algebras. Let us start with the following general situation.

Let $\Phi$ be a subset of ${\mathbb N}$. Given an ${\mathbb
N}$-graded $R$-algebra $\Lambda=\bigoplus_{i\ge
0}^{\infty}\Lambda_i$, where $R$ is a commutative ring and each
$\Lambda_i$ is an $R$-module with $\Lambda_i\Lambda_j\subseteq
\Lambda_{i+j}$ for all $i,j\in {\mathbb N}$, we define an $R$-module
$\Lambda(\Phi):=\bigoplus_{i\in \Phi} \Lambda_i$, and a
multiplication in $A(\Phi)$: for $a_i\in \Lambda_i$ and $b_j\in
\Lambda_j$ with $i,j\in \Phi$, we define $a_i\cdot b_j = a_ib_j$ if
$i+j\in \Phi$, and zero otherwise. Then one can easily check that
$\Lambda(\Phi)$ is an associative algebra if $\Phi$ is an admissible
subset of ${\mathbb N}$.

This procedure can be applied to a triangulated category, in this
special situation, the details which are needed in our proofs read
as follows:

Let $\mathcal T$ be a triangulated $R$-category over a commutative
Artin ring $R$, and let $\Phi$ be a subset in $\mathbb{N}$
containing $0$. We denote by $E_{\mathcal T}^{\Phi}(-,-)$ the
bifunctor
$$\bigoplus_{i\in\Phi}\, \Hom_{\cal T}(-, -[i]): {\mathcal T}\times{\mathcal T}\lra R\mbox{\rm -Mod},$$
$$ (X,Y)\mapsto E_{\mathcal T}^{\Phi}(X,Y):= \bigoplus_{i\in\Phi}\, \Hom_{\mathcal T}(X, Y[i]).$$
Let $X, Y$ and $Z$ be objects in $\mathcal T$. For each $i\in\Phi$,
let $\iota_i$ denote the canonical embedding of $\Hom_{\cal
T}(X,Y[i])$ into $E_{\mathcal T}^{\Phi}(X,Y).$ For $i\not\in\Phi$,
we define $\iota_i$ to be the zero map from $\Hom_{\cal T}(X,Y[i])$
to $E_{\cal T}^{\Phi}(X,Y)$. An element in $E_{\cal T}^{\Phi}(X,Y)$
is of the form $(f_i)_{i\in\Phi}$, where $f_i$ is a morphism in
$\Hom_{\cal T}(X,Y[i])$ for $i\in\Phi$. For simplicity, we shall
just write $(f_i)$ for $(f_i)_{i\in\Phi}$, and each element $(f_i)$
in $E_{\cal T}^{\Phi}(X,Y)$ can be rewritten as
$\displaystyle{\sum_{i\in\Phi}}\iota_i(f_i)$, where $\iota_i(f_i)$
denotes the image of $f_i$ under the map $\iota_i$.

Let $(f_i)\in\mbox{E}_{\cal T}^{\Phi}(X,Y)$ and
$(g_i)\in\mbox{E}_{\cal T}^{\Phi}(Y,Z)$. We define a multiplication
$(h_i)=(f_i)(g_i)$:

$$\mbox{E}_{\cal T}^{\Phi}(X,Y)\times\mbox{E}_{\cal T}^{\Phi}(Y,Z)\lra\mbox{E}_{\cal T}^{\Phi}(X,Z)$$
$$ \big((f_i),(g_i)\big)\mapsto (h_i),$$ where
$$h_i:=\sum_{{{u, v\in\Phi}\atop {u+v=i}}}f_u(g_v[u])$$
for each $i\in\Phi$. In particular, for $f\in\Hom_{\cal T}(X,Y[i])$
and $g\in\Hom_{\cal T}(X,Y[j])$ with $i,j\in \Phi$, we have
$$\iota_i(f) \; \iota_j(g)=\iota_{i+j}(f(g[i])).$$
Note that $\iota_{i+j}=0$ if $i+j\not\in\Phi$.

The next proposition explains further why we need to introduce
admissible subsets.

\begin{Prop}
   Let ${\cal T}$ be a triangulated $R$-category with at least one non-zero object, and let $\Phi$ be a
subset of $\mathbb{N}$ containing $0$. Then $\Ex{{\cal T}}{\Phi}{V}$
together with the multiplication defined above is an associative
$R$-algebra for every object $V\in{\cal T}$ if and only if $\Phi$ is
an admissible subset of $\mathbb N$.
\end{Prop}

{\it Proof.} If $\Phi$ is an admissible subset of $\mathbb{N}$, then
it is straightforward to check that the multiplication on $\Ex{{\cal
T}}{\Phi}{V}$ defined above is associative for all objects
$V\in{\cal T}$. Now we assume that $\Phi$ is not an admissible
subset, that is, there are integers $i,j, k\in\Phi$ satisfying:
 $i+j+k\in\Phi$,  $i+j\in\Phi$, and $j+k\not\in\Phi$. Let $X$ be a
non-zero object in ${\cal T}$, and let
$V:=\bigoplus_{s=0}^{i+j+k}X[s]$. We consider the multiplication on
$\Ex{\cal T}{\Phi}{V}$. By definition, the object
$\bigoplus_{s=i}^{i+j+k}X[s]$ is a common direct summand of $V$ and
$V[i]$. Let $f$ be the composition $V\lraf{\pi}
\bigoplus_{s=i}^{i+j+k}X[s]\lraf{\lambda} V[i]$, where $\pi$ is the
canonical projection and $\lambda$ is the canonical inclusion.
Similarly, we define $g: V\lra \bigoplus_{s=j}^{i+j+k}X[s]\lra V[j]$
and $h: V\lra \bigoplus_{s=k}^{i+j+k}X[s]\lra V[k]$. Since
$i+j\in\Phi$, we have
$\big(\iota_i(f)\iota_j(g)\big)\iota_k(h)=\iota_{i+j}(f(g[i]))\iota_k(h)=\iota_{i+j+k}\big(f(g[i])(h[i+j])\big)$.
One can check that $f(g[i])(h[i+j])$ is just the composition $V\lra
X[i+j+k]\lra V[i+j+k]$, where the maps are canonical maps. Hence the
map $\big(\iota_i(f)\iota_j(g)\big)\iota_k(h)$ is non-zero. Since
$j+k\not\in\Phi$, we have $\iota_j(g)\iota_k(h)=0$, and consequently
$\iota_i(f)\big(\iota_j(g)\iota_k(h)\big)=0$. This shows that the
multiplication of $\Ex{\cal T}{\Phi}{V}$ is not associative, and the
proof is completed. $\square$

Note that $E^{{\mathbb N}}_{\mathcal T}(X)$ is an ${\mathbb
N}$-graded associative $R$-algebra with $\Hom_{\mathcal T}(X,X[i])$
as $i$-th component. If we define $\Lambda := E^{{\mathbb
N}}_{\mathcal T}(X)$, then $\Lambda(\Phi) = E^{\Phi}_{\mathcal
T}(X)$.

From now on, we consider exclusively admissible subsets $\Phi$ of
$\mathbb N$. Thus, for objects ${X}$ and $Y$ in $\mathcal T$, one
has an $R$-algebra $\mbox{E}_{\cal T}^{\Phi}(X,X)$ (which may not be
artinian), and a left $\mbox{E}_{\cal T}^{\Phi}(X,X)$-module
$\mbox{E}_{\cal T}^{\Phi}(X,Y)$. For simplicity, we write
$\mbox{E}_{\cal T}^{\Phi}(X)$ for $\mbox{E}_{\cal T}^{\Phi}(X,X)$.

In case $\Phi=\Phi(1,0)$, we see that $\mbox{E}_{\cal T}^{\Phi}(X)$
is the endomorphism algebra of the object $X$ in $\mathcal T$.  In
case $\Phi={\mathbb N}$, we know that $\mbox{E}_{\cal T}^{\Phi}(X)$
is the generalized Yoneda algebra Ext$^*_{\cal T}(X)=\bigoplus_{i\ge
0}\Hom_{\cal T}(X,X[i])$ of $X$. Particularly, let us take $\cal T$
= $\Db{A}$ with $A$ an Artin $R$-algebra. If $A$ is
representation-finite and if $X$ is an additive generator for
$A$-mod, then $\Ex{\cal T}{\Phi(1,0)}{X}$ is the Auslander algebra
of $A$; if we take $X=A/\rad(A)$, then
$\mbox{E}^{\Phi(1,\infty)}_{\mathcal T}(X)$ is the usual Yoneda
algebra of $A$. Thus the algebra $\mbox{E}_{\cal T}^{\Phi}(X)$ is a
generalization of both Auslander algebra and Yoneda algebra. For
this reason, the algebra $\mbox{E}_{\cal T}^{\Phi}(X)$ of $X$ in a
triangulated category $\mathcal T$ is called the
$\Phi$-\emph{Auslander-Yoneda algebra} of $X$ in $\mathcal T$ in
this paper.

If $\mathcal T$ = $\Db{A}$ with $A$ an Artin algebra, we simply
write $\mbox{E}_{A}^{\Phi}(X)$ for $\mbox{E}_{\cal T}^{\Phi}(X)$,
and $\mbox{E}_{A}^{\Phi}(X,Y)$ for $\mbox{E}_{\cal T}^{\Phi}(X,Y)$.
If $\Phi$ is finite, or if the projective or injective dimension of
$X$ is finite, then $\Ex{A}{\Phi}{X}$ is an Artin $R$-algebra.

Note also that the algebra $\mbox{E}_{\cal T}^{\Phi(1,m)}(X)$ is a
quotient algebra of $\mbox{E}_{\cal T}^{\mathbb N}(X)$, and the
algebra $\mbox{E}_{\cal T}^{\Phi(n, m)}(X)$ is a subalgebra of
$\mbox{E}_{\cal T}^{\Phi(1, nm)}(X)$. Nevertheless, if we take
$\Phi=\{0,3,9\}$ and $X$ a simple module over the algebra
$A:=k[X]/(X^2)$ with $k$ a field, then $\Ex{A}{\Phi}{X}$ is neither
a subalgebra nor a quotient algebra of the generalized Yoneda
algebra of $X$.

Let us remark that one may define the notion of an admissible subset
of $\mathbb Z$ (or of a monoid $M$ with an identity $e$), and
introduce $\Phi$-Auslander-Yoneda algebra of an object in an
arbitrary $R$-category $\cal C$ with an additive self-equivalence
functor (or a family of additive functors $\{F_g\}_{g\in M}$ from
$\cal C$ to itself, such that $F_e=id_{\mathcal C}$ and
$F_gF_h=F_{gh}$ for all $g,h\in M$). For our goals in this paper, we
just formulate the admissible subsets for $\mathbb N$.

\subsection{Almost $\nu$-stable derived equivalences}

We briefly recall some basic facts on almost $\nu$-stable derived
equivalences from \cite{HX3}, which are needed in proofs.

\medskip
Let $A$ and $B$ be Artin algebras, and let $F:\Db{A}\lra\Db{B}$ be a
derived equivalence between $A$ and $B$. Suppose that $\cpx{Q}$ and
$\cpx{\bar{Q}}$ are the tilting complexes associated to $F$ and to a
quasi-inverse $G$ of $F$, respectively. Now, we assume that $Q^i=0$
for all $i>0$, that is, the complex $\cpx{Q}$ is of the form
$$0\lra Q^{-n}\lra \cdots\lra Q^{-1}\lra Q^0\lra 0.$$
In this case, the complex $\cpx{\bar{Q}}$ may be chosen of the
following form (see \cite[Lemma 2.1]{HX3}, for example)
$$0\lra \bar{Q}^{0}\lra \bar{Q}^{1}\lra\cdots\lra \bar{Q}^{n}\lra 0.$$
Set $Q:=\bigoplus_{i=1}^nQ^{-i}$ and
$\bar{Q}:=\bigoplus_{i=1}^n\bar{Q}^i$. The functor $F$ is called an
\emph{almost $\nu$-stable} derived equivalence provided
$\add(_AQ)=\add(\nu_AQ)$ and $\add(_B\bar{Q})=\add(\nu_B\bar{Q})$. A
crucial property is that an almost $\nu$-stable derived equivalence
induces an equivalence between the stable module categories
$\stmodcat{A}$ and $\stmodcat{B}$. Thus $A$ and $B$ share many
common properties, for example, $A$ is representation-finite if and
only if $B$ is representation-finite.

In the following lemma, we collect some basic facts on almost
$\nu$-stable derived equivalences, which will be used in our proofs.

\begin{Lem}\label{AlmostV-property}
 Let $F:\Db{A}\ra\Db{B}$ be an almost $\nu$-stable derived
 equivalence between Artin algebras $A$ and $B$. Suppose that
 $\cpx{Q}$ and $\cpx{\bar{Q}}$ are the tilting complexes associated to
 $F$ and to its quasi-inverse $G$, respectively. Then:

$(1)$ For each $A$-module $X$, the complex $F(X)$ is isomorphic in
$\Db{B}$ to a radical complex $\cpx{\bar{Q}_X}$ of the form
$$0\lra \bar{Q}_X^0\lra \bar{Q}_X^1\lra\cdots\lra \bar{Q}_X^n\lra 0,$$
with $\bar{Q}_X^i\in\add(_B\bar{Q})$ for all $i>0$. Moreover, the
complex $\cpx{\bar{Q}_X}$ of this form is unique up to isomorphism
in $\Cb{B}$.

$(2)$ For each $B$-module $Y$, the complex $G(Y)$ is isomorphic in
$\Db{A}$ to a radical complex $\cpx{Q_Y}$ of the form
$$0\lra {Q}_Y^{-n}\lra\cdots\lra {Q}_Y^{-1}\lra {Q}_Y^0\lra 0,$$
with $Q_Y^i\in\add(_AQ)$ for all $i<0$. Moreover, the complex
$\cpx{Q_Y}$ of this form is unique up to isomorphism in $\Cb{B}$.

$(3)$ There is a stable equivalence
$\bar{F}:\stmodcat{A}\lra\stmodcat{B}$ with $\bar{F}(X)=\bar{Q}_X^0$
for each $A$-module $X$.

$(4)$ There is a stable equivalence $\bar{G}:
\stmodcat{B}\lra\stmodcat{A}$ with $\bar{G}(Y)=Q_Y^0$ for each
$B$-module $Y$. Moreover, the functor $\bar{G}$ is a quasi-inverse
of $\bar{F}$ defined in $(3)$.

$(5)$ For an $A$-module $X$,  we denote by $\bar{Q}_X^+$  the
complex $\sigma_{>0}\cpx{\bar{Q}_X}$. Then $G(\bar{Q}_X^+)$ is
isomorphic in $\Db{A}$ to a bounded complex $\cpx{P_X}$ of
projective-injective $A$-modules with $P_X^i=0$ for all $i>1$.
\end{Lem}

{\it Proof.} The statement (1) follows from \cite[Lemma 3.1]{HX3}.
The statement (2) is a direct consequence of the definition of an
almost $\nu$-stable derived equivalence and \cite[Lemma 3.2]{HX3}.
Note that the statements (3) and (4) follow from the proof of
\cite[Theorem 3.7]{HX3}, and the statement (5) is implied in the
proof of \cite[Proposition 3.6]{HX3}. $\square$

\medskip
For an Artin algebra $A$, let $_{A}E$ be the direct sum of all
non-isomorphic indecomposable projective $A$-modules $P$ with the
property: $\nu_{A}^iP$ is projective-injective for all $i\geq 0$.
The $A$-module $_{A}E$ is called a \emph{maximal}
$\nu$-\emph{stable} $A$-module.

\subsection{Derived equivalences for Auslander-Yoneda algebras}

\medskip
Our main result in this section is the following theorem on derived
equivalences between $\Phi$-Auslander-Yoneda algebras.
\begin{Theo}
Let $F:\Db{A}\lra\Db{B}$ be an almost $\nu$-stable derived
equivalence between two Artin algebras $A$ and $B$, and let
$\bar{F}$ be the stable equivalence defined in {\rm Lemma
\ref{AlmostV-property} (3)}. For an $A$-module $X$, we set
$M:=A\oplus X$ and $N:=B\oplus \bar{F}(X)$. Suppose that $\Phi$ is
an admissible subset in $\mathbb{N}$. Then we have the following:

\smallskip
$(1)$ The algebras $\Ex{A}{\Phi}{M}$ and $\Ex{B}{\Phi}{N}$ are
derived-equivalent.

$(2)$ If $\Phi$ is finite, then there is an almost $\nu$-stable
derived equivalence between $\Ex{A}{\Phi}{M}$ and $\Ex{B}{\Phi}{N}$.
Thus $\Ex{A}{\Phi}{M}$ and $\Ex{B}{\Phi}{N}$ are also stably
equivalent. In particular, there is an almost $\nu$-stable derived
equivalence and a stable equivalence between $\End_A(M)$ and
$\End_B(N)$.\label{ThmYonedaAlgebra}
\end{Theo}

Thus, under the assumptions of Theorem \ref{ThmYonedaAlgebra}, if
$\Phi$ is finite, then the algebras $\Ex{A}{\Phi}{M}$ and
$\Ex{B}{\Phi}{N}$ share many common invariants: for example,
finiteness of finitistic and global dimensions, representation
dimension, Hochschild cohomology, representation-finite type and so
on.

\medskip The rest of this section is
essentially devoted to the proof of Theorem \ref{ThmYonedaAlgebra}.
First of all, we need some preparations. Let us start with the
following lemma that describes some basic properties of the algebra
$\Ex{A}{\Phi}{V}$, where $V$ is an $A$-module and is considered as a
complex concentrated on degree zero.

\begin{Lem}
Let $A$ be an Artin algebra, and let $V$ be an $A$-module. Suppose
that $V_1$ and $V_2$ are in $\add(_AV)$. Then

\medskip
$(1)$ The $\Ex{A}{\Phi}{V}$-module $\E{A}{\Phi}{V}{V_1}$ is
projective and finitely generated, and there is an isomorphism
$$\mu: \E{A}{\Phi}{V_1}{V_2}\lra \Hom_{\Ex{A}{\Phi}{V}}(\E{A}{\Phi}{V}{V_1}, \E{A}{\Phi}{V}{V_2}),$$
which sends $(f_i)\in\E{A}{\Phi}{V_1}{V_2}$ to the morphism
$\big((a_i)\mapsto (a_i)(f_i)\big)$ for $(a_i)\in
\E{A}{\Phi}{V}{V_2}$. Moreover, if $V_3\in\add(_AV)$ and
$(g_i)\in\E{A}{\Phi}{V_2}{V_3}$, then
$\mu((f_i)(g_i))=\mu((f_i))\mu((g_i))$.

\smallskip
$(2)$ The functor $\E{A}{\Phi}{V}{-}: \add(_AV)\lra
\pmodcat{\Ex{A}{\Phi}{V}}$ is faithful.

\smallskip
$(3)$ If $V_1$ is projective or $V_2$ is injective, then the functor
$\E{A}{\Phi}{V}{-}$ induces an isomorphism of $R$-modules:
$$\E{A}{\Phi}{V}{-}: \Hom_A(V_1, V_2)\lra \Hom_{\Ex{A}{\Phi}{V}}(\E{A}{\Phi}{V}{V_1}, \E{A}{\Phi}{V}{V_2}).$$

$(4)$ If $\Phi$ is finite, and  $P\in\add(_AV)$ is projective, then
$$\nu_{\Ex{A}{\Phi}{V}}\E{A}{\Phi}{V}{P}\simeq \E{A}{\Phi}{V}{\nu_AP}.$$
   \label{Extalgproperty}
\end{Lem}
{\it Proof.}
 (1) Since $\E{A}{\Phi}{V}{-}$ is an additive functor and since
 $V_1\in\add(_AV)$, we know that $\E{A}{\Phi}{V}{V_1}$ is in
 $\add(\Ex{A}{\Phi}{V})$, and consequently $\E{A}{\Phi}{V}{V_1}$ is a finitely generated projective
 $\Ex{A}{\Phi}{V}$-module. Similarly, the $\Ex{A}{\Phi}{V}$-module
 $\E{A}{\Phi}{V}{V_2}$ is also projective. To show that $\mu$ is an
 isomorphism,  we can assume that $V_1$ is
 indecomposable by additivity. Let $\pi_1: V\lra V_1$ be the canonical projection,
 and let $\lambda_1: V_1\lra V$ be the canonical injection. We define a map
 $$\gamma: \Hom_{\Ex{A}{\Phi}{V}}(\E{A}{\Phi}{V}{V_1}, \E{A}{\Phi}{V}{V_2})\lra \E{A}{\Phi}{V_1}{V_2}$$
 by sending $\alpha\in \Hom_{\Ex{A}{\Phi}{V}}(\E{A}{\Phi}{V}{V_1}, \E{A}{\Phi}{V}{V_2})$ to $\iota_0(\lambda_1)\alpha\big(\iota_0(\pi_1)\big)$.
By calculation, the morphism $(\gamma\mu)(\alpha):
\E{A}{\Phi}{V}{V_1}\lra \E{A}{\Phi}{V}{V_2}$ sends each
$x\in\E{A}{\Phi}{V}{V_1}$ to $x\iota_0(\lambda_1)\alpha
\big(\iota_0(\pi_1)\big)=\alpha\big(x\iota_0(\lambda_1)\iota_0(\pi_1)\big)=\alpha(x)$.
This shows that $\gamma\,\mu={\rm id}$. Similarly, one can check
that $\mu\,\gamma={\rm id}$. Hence $\mu$ is an isomorphism. The rest
of (1) can be verified easily.

(2) Using definition, one can check that the map
$$\E{A}{\Phi}{V}{-}:
\Hom_{\Db{A}}(V_1, V_2)\lra
\Hom_{\Ex{A}{\Phi}{V}}(\E{A}{\Phi}{V}{V_1}, \E{A}{\Phi}{V}{V_2})$$
is the composition of the embedding $\iota_0: \Hom_A(V_1, V_2)\lra
\E{A}{\Phi}{V_1}{V_2}$ with the isomorphism $\mu$ in (1). Hence
$\E{A}{\Phi}{V}{-}$ is a faithful functor.

(3) If $V_1$ is projective or $V_2$ is injective, then the embedding
$$\iota_0: \Hom_{\Db{A}}(V_1, V_2)\lra \E{A}{\Phi}{V_1}{V_2}$$
is an isomorphism. Since $\E{A}{\Phi}{V}{-}$ is the composition of
$\iota_0$ with the isomorphism $\mu$ in (1), the statement (3)
follows.

(4) This follows directly from the following isomorphisms
$$\begin{array}{rcl}
\nu_{\Ex{A}{\Phi}{V}}\E{A}{\Phi}{V}{P} & = &
D\Hom_{\Ex{A}{\Phi}{V}}(\E{A}{\Phi}{V}{P},
\E{A}{\Phi}{V}{V})\\
& \simeq & D\E{A}{\Phi}{P}{V}  \quad \mbox{by (1)}\\
& = & D\Hom_A(P, V)\\
&\simeq & \Hom_A(V, \nu_AP)\\
& =& \E{A}{\Phi}{V}{\nu_AP}.
\end{array}$$
Thus we have finished the proof. $\square$

From now on, we assume that  $F: \Db{A}\lra \Db{B}$ is an almost
$\nu$-stable derived equivalence with a quasi-inverse functor $G$,
that $\cpx{Q}$ and $\cpx{\bar{Q}}$ are tilting complexes associated
to $F$ and $G$, respectively, and that $\bar{F}:
\stmodcat{A}\lra\stmodcat{B}$ is the stable equivalence defined by
Lemma \ref{AlmostV-property} (3). For an $A$-module $X$, we may
assume that $F(X)=\cpx{\bar{Q}_X}$ as in Lemma
\ref{AlmostV-property} (1), and define $_AM=A\oplus X$ and
$_BN=B\oplus \bar{F}(X)$. By $\cpx{\bar{T}}$ we denote the complex
$\cpx{\bar{Q}}\oplus\cpx{\bar{Q}_X}$. Clearly, $\cpx{\bar{T}}$ is in
$\Kb{\add(_BN)}$.

\begin{Lem}
Keeping the notations above, we have the following:

$(1)$ $\Hom_{\Kb{\add(_BN)}}(\cpx{\bar{T}}, \cpx{\bar{T}}[i])=0$ for
all $i\neq 0$.

$(2)$ $\add(\cpx{\bar{T}})$ generates $\Kb{\add(_BN)}$ as a
triangulated category. \label{barT-orth}
\end{Lem}

{\it Proof.} Since $F(A)\simeq \cpx{\bar{Q}}$, the complex
$\cpx{\bar{T}}$ is isomorphic to $F(M)=\cpx{\bar{Q}_{M}}$. So, we
consider $\cpx{\bar{Q}_M}$ instead.

\smallskip
(1)  Suppose $i<0$. Then $\Hom_{\Kb{B}}(\cpx{\bar{Q}_M},
\cpx{\bar{Q}_M}[i])\simeq\Hom_{\Db{B}}(\cpx{\bar{Q}_M},
\cpx{\bar{Q}_M}[i])$ by Lemma \ref{kdiso}. Since
$F(M)=\cpx{\bar{Q}_M}$, we have $\Hom_{\Db{B}}(\cpx{\bar{Q}_M},
\cpx{\bar{Q}_M}[i])\simeq\Hom_{\Db{A}}(M, M[i])=0$. Hence
$\Hom_{\Kb{B}}(\cpx{\bar{Q}_M}, \cpx{\bar{Q}_M}[i])$ = $0$ for all
$i<0$.

Let $\bar{Q}_M^+$ be the complex $\sigma_{>0}\cpx{\bar{Q}_M}$. There
is a distinguished triangle
$$(*)\hspace{.5cm} \bar{Q}_M^+\stackrel{i_M}{\lra}
\cpx{\bar{Q}_M}\stackrel{\pi_M}{\lra}
\bar{F}(M)\stackrel{\alpha_M}{\lra}\bar{Q}_M^+[1]$$ in $\Kb{B}$.
%, where $\alpha_M$ is given by the differential map from
%$\bar{Q}_M^0$ to $\bar{Q}_M^1$.
Applying $\Hom_{\Kb{B}}(\cpx{\bar{Q}_M},- )$ to $(*)$, we get an
exact sequence
 {\small $$ \Hom_{\Kb{B}}(\cpx{\bar{Q}_M},
\bar{F}(M)[i-\!1])\ra\Hom_{\Kb{B}}(\cpx{\bar{Q}_M},
\bar{Q}_M^+[i])\ra \Hom_{\Kb{B}}(\cpx{\bar{Q}_M},
\cpx{\bar{Q}_M}[i])\ra\Hom_{\Kb{B}}(\cpx{\bar{Q}_M},
\bar{F}(M)[i])$$}
 for each integer $i$. Since $\bar{Q}_M^i=0$ for all
$i< 0$, we have $\Hom_{\Kb{B}}(\cpx{\bar{Q}_M}, \bar{F}(M)[i])=0$
for all $i>0$. By Lemma \ref{AlmostV-property} (5), $G(\bar{Q}_M^+)$
is isomorphic to a bounded complex $\cpx{P_M}$ of
projective-injective $A$-modules such that $P_M^i=0$ for all $i>1$.
Thus, we have
$$\begin{array}{rcl}
\Hom_{\Kb{B}}(\cpx{\bar{Q}_M}, \bar{Q}_M^+[i]) & \simeq &
\Hom_{\Db{B}}(\cpx{\bar{Q}_M}, \bar{Q}_M^+[i]) \quad\quad (\mbox{since\;} \bar{Q}_M^+[i] \mbox{ is in }\Kb{\imodcat{B}})\\
& \simeq &
\Hom_{\Db{A}}(G(\cpx{\bar{Q}_M}), G(\bar{Q}_M^+)[i])\\
&\simeq & \Hom_{\Db{A}}(M, \cpx{P_M}[i])\\
&\simeq & \Hom_{\Kb{A}}(M, \cpx{P_M}[i]) \\ & = & 0
\end{array}$$
for all $i>1$, and consequently
$\Hom_{\Kb{B}}(\cpx{\bar{Q}_M},\cpx{\bar{Q}_M}[i])=0$ for all $i>1$.
To prove (1), it remains to show that
$\Hom_{\Kb{B}}(\cpx{\bar{Q}_M}, \cpx{\bar{Q}_M}[1])=0$. Using the
above exact sequence, we only need to show that the induced map
$$\Hom_{\Kb{B}}(\cpx{\bar{Q}_M},\alpha_M):
\Hom_{\Kb{B}}(\cpx{\bar{Q}_M},
\bar{F}(M))\lra\Hom_{\Kb{B}}(\cpx{\bar{Q}_M}, \bar{Q}_M^+[1])$$
 is surjective. Note that $G(\bar{Q}_M^+)$ is
isomorphic in $\Db{A}$ to a complex $\cpx{P_M}$ of
projective-injective modules such that $P_M^k=0$ for all $k>1$.
Thus, we can form a commutative diagram
$$\begin{CD}
\cpx{P_M}@>{\phi_M}>> M@>{\lambda}>>\con(\phi_M)@>p>>\cpx{P_M}[1]\\
@VV{\simeq}V  @VV{\simeq}V  @VV{\simeq}V  @VV{\simeq}V\\
G({\bar{Q}_M^+})@>{Gi_M}>>GF(M)@>{G\pi_M}>>G\bar{F}(M)@>G\alpha_M>> G({\bar{Q}_M^+})[1]\\
\end{CD}$$
in $\Db{A}$, where the vertical maps are all isomorphisms, $\lambda$
and $p$ are the canonical morphisms, and where the morphism $\phi_M$
is chosen in $\Kb{A}$ such that the first square is commutative. The
distinguished triangle in the top row of the above diagram can be
viewed as a distinguished triangle in $\Kb{A}$. Applying
$\Hom_{\Kb{A}}(M, -)$ to this triangle, we can easily see that
$\Hom_{\Kb{A}}(M, p)$ is a surjective map since $\Hom_{\Kb{A}}(M,
M[1])=0$. By Lemma 2.1, the localization functor $\theta:
\Kb{A}\ra\Db{A}$ induces two isomorphisms
$$\Hom_{\Kb{A}}(M, \con(\phi_M))\simeq \Hom_{\Db{A}}(M,
\con(\phi_M))\; \mbox{and}\; \Hom_{\Kb{A}}(M, \cpx{P_M}[1])\simeq
\Hom_{\Db{A}}(M, \cpx{P_M}[1]).$$ It follows that $\Hom_{\Db{A}}(M,
p)$ is surjective. Since the vertical maps of the above diagram are
all isomorphisms, the map $\Hom_{\Db{A}}(M, G\alpha_M)$ is
surjective, or equivalently $\Hom_{\Db{A}}(G(\cpx{\bar{Q}}_M)$,
 $G\alpha_M)$ is surjective. Since $G$ is an equivalence, it follows
that $\Hom_{\Db{B}}(\cpx{\bar{Q}}_M, \alpha_M)$ is surjective. By
Lemma 2.1 again, the localization functor $\theta: \Kb{B}\ra \Db{B}$
gives rise to isomorphisms $$\Hom_{\Kb{B}}(\cpx{\bar{Q}_M},
\bar{F}(M))\simeq \Hom_{\Db{B}}(\cpx{\bar{Q}_M},
\bar{F}(M))\;\mbox{and} \; \Hom_{\Kb{B}}(\cpx{\bar{Q}_M},
\bar{Q}_M^+[1])\simeq \Hom_{\Db{B}}(\cpx{\bar{Q}_M},
\bar{Q}_M^+[1]).$$ Hence the map $\Hom_{\Kb{B}}(\cpx{\bar{Q}}_M,
\alpha_M)$ is surjective, and consequently
$\Hom_{\Kb{B}}(\cpx{\bar{Q}_M}, \cpx{\bar{Q}_M}[1])=0$.

Altogether, we have shown that $\Hom_{\Kb{B}}(\cpx{\bar{Q}_M},
\cpx{\bar{Q}_M}[i])=0$ for all $i\neq 0$. Since $\Kb{B}$ is a full
subcategory of $\Kb{\add(_BN)}$, we have
$\Hom_{\Kb{\add(_BN)}}(\cpx{\bar{Q}_M}, \cpx{\bar{Q}_M}[i])=0$ for
all $i\neq 0$. This proves (1).

\smallskip
(2) Since $\cpx{\bar{Q}}$ is a tilting complex over $B$,
$\add(\cpx{\bar{Q}})$ generates $\Kb{\add(_BB)}$ as a triangulated
category. By Lemma \ref{AlmostV-property}, $\bar{Q}_X^0=\bar{F}(X)$
and all the terms of $\cpx{\bar{Q}_X}$ other than $\bar{Q}_X^0$ are
in $\add(_BB)$. Hence $\bar{F}(X)$ is in the triangulated
subcategory generated by $\add(\cpx{\bar{Q}}\oplus\cpx{\bar{Q}_X})$,
and consequently $\add(\cpx{\bar{Q}}\oplus\cpx{\bar{Q}_X})$
generates $\Kb{\add(B\oplus\bar{F}(X))}$ as a triangulated category.
Thus, the statement (2) follows. $\square$

\medskip
 The additive functor $\E{B}{\Phi}{N}{-}: \add({}_BN)\lra
\pmodcat{\Ex{B}{\Phi}{N}}$ induces a triangle functor
$$\EP{B}{\Phi}{N}{-}: \Kb{\add({}_BN)}\lra \Kb{\pmodcat{\Ex{B}{\Phi}{N}}}.$$
For each integer $i$, the $i$-th term of
$\EP{B}{\Phi}{N}{\cpx{\bar{T}}}$ is $\E{B}{\Phi}{N}{\bar{T}^i}$, and
the differential map from $\E{B}{\Phi}{N}{\bar{T}^i}$ to
$\E{B}{\Phi}{N}{\bar{T}^{i+1}}$ is $\E{B}{\Phi}{N}{d}$, where $d:
\bar{T}^i\lra \bar{T}^{i+1}$ is the differential map of
$\cpx{\bar{T}}$.

\medskip
\begin{Lem}
  The complex $\EP{B}{\Phi}{N}{\cpx{\bar{T}}}$ is a tilting complex over
  $\Ex{B}{\Phi}{N}$.\label{lemtiltYoneda}
\end{Lem}

\medskip
{\it Proof.}
 Let $i\neq 0$, and let $\cpx{f}$ be a morphism in
$\Hom_{\Kb{\pmodcat{\Ex{B}{\Phi}{N}}}}(\EP{B}{\Phi}{N}{\cpx{\bar{T}}},
 \EP{B}{\Phi}{N}{\cpx{\bar{T}}}[i])$. Then we have a commutative diagram
$$\xymatrix{
 0 \ar[r] &
   \E{B}{\Phi}{N}{\bar{T}^{0}}\ar[r]^{\E{B}{\Phi}{N}{d}}\ar[d]^{f^0} &
   \E{B}{\Phi}{N}{\bar{T}^{1}}\ar[r]^(.6){\E{B}{\Phi}{N}{d}}\ar[d]^{f^1} & \cdots\\
 \E{B}{\Phi}{N}{\bar{T}^{i-1}}\ar[r]^{\E{B}{\Phi}{N}{d}} &
   \E{B}{\Phi}{N}{\bar{T}^{i}}\ar[r]^{\E{B}{\Phi}{N}{d}} &
   \E{B}{\Phi}{N}{\bar{T}^{i+1}}\ar[r]^(.6){\E{B}{\Phi}{N}{d}} & \cdots .\\
  }$$
Note that the term $\E{B}{\Phi}{N}{\bar{T}^{i}}$ is zero if $i<0$.
Since all the terms of $\cpx{\bar{T}}$ other than $\bar{T}^0$ are
projective-injective, and since $i\neq 0$,  we see from Lemma
\ref{Extalgproperty} (3) that $f^k=\E{B}{\Phi}{N}{g^k}$ for some
$g^k: \bar{T}^k\lra\bar{T}^{k+i}$ for all integers $k$. It follows
from the above commutative diagram that, for each integer $k$, we
have
$$\E{B}{\Phi}{N}{d}\E{B}{\Phi}{N}{g^{k+1}}-\E{B}{\Phi}{N}{g^k}\E{B}{\Phi}{N}{d}=0, $$
or equivalently $\E{B}{\Phi}{N}{dg^{k+1}-g^kd}=0$.  Since
$\E{B}{\Phi}{N}{-}: \add(_BN)\lra\pmodcat{\Ex{B}{\Phi}{N}}$ is a
faithful functor by Lemma \ref{Extalgproperty} (2), we have
$dg^{k+1}-g^kd=0$ for all integers $k$, and consequently
$\cpx{g}:=(g^k)$ is in $\Hom_{\Kb{\add(_BN)}}(\cpx{\bar{T}},
\cpx{\bar{T}}[i])$ and $\cpx{f}=\EP{B}{\Phi}{N}{\cpx{g}}$. By Lemma
\ref{barT-orth} (1), the map $\cpx{g}$ is null-homotopic, and
consequently $\cpx{f}=\EP{B}{\Phi}{N}{\cpx{g}}$ is null-homotopic.
Thus, we have proved that
$$\Hom_{\Kb{\pmodcat{\Ex{B}{\Phi}{N}}}}(\EP{B}{\Phi}{N}{\cpx{\bar{T}}},
 \EP{B}{\Phi}{N}{\cpx{\bar{T}}}[i])=0$$
for all non-zero integers $i$.

\medskip
By definition, the triangle functor $\EP{B}{\Phi}{N}{-}:
\Kb{\add(_BN)}\lra\Kb{\pmodcat{\Ex{B}{\Phi}{N}}}$ sends $N$ to
$\Ex{B}{\Phi}{N}$. The full triangulated subcategory of
$\Kb{\add(_BN)}$ generated by $\add(\cpx{\bar{T}})$ contains $N$ by
Lemma \ref{barT-orth} (2), and so $\Ex{B}{\Phi}{N}$ is in the full
triangulated subcategory of $\Kb{\pmodcat{\Ex{B}{\Phi}{N}}}$
generated by $\add(\EP{B}{\Phi}{N}{\cpx{\bar{T}}})$. Hence
$\add(\EP{B}{\Phi}{N}{\cpx{\bar{T}}})$ generates
$\Kb{\pmodcat{\Ex{B}{\Phi}{N}}}$ as a triangulated category. This
finishes the proof. $\square$

In the following, we shall prove that the endomorphism algebra of
the complex $\EP{B}{\Phi}{N}{\cpx{\bar{T}}}$ is isomorphic to
$\Ex{A}{\Phi}{M}$. For this purpose, we first prove the following
lemma.

\begin{Lem} Keeping the notations above, for each $A$-module $V$, we have:

$(1)$ For each positive integer $k$, there is an isomorphism
$$\theta_k: \Hom_{\Db{A}}(V, V[k])\lra \Hom_{\Db{B}}(\bar{F}(V),
\bar{F}(V)[k]).$$ Here we denote the image of $g$ under $\theta_k$
by $\theta_k(g)$.

$(2)$ For each pair of positive integers $k$ and $l$, the $\theta_k$
and $\theta_l$ in $(1)$ satisfy
$$\theta_k(g)(\theta_l(h)[k])=\theta_{k+l}(g(h[k]))$$
for all $g\in\Hom_{\Db{A}}(V, V[k])$ and $h\in\Hom_{\Db{A}}(V,
V[l])$. \label{thetamaps}
\end{Lem}

\medskip
{\it Proof.} By Lemma \ref{AlmostV-property}, we may assume that
$F(V)$ is the complex $\cpx{\bar{Q}_V}$ defined in Lemma
\ref{AlmostV-property} (1), and therefore $\bar{F}(V)=\bar{Q}_V^0$.
As before, the complex $\sigma_{>0}\cpx{\bar{Q}_V}$ is denoted by
$\bar{Q}_V^+$. Thus, we have a distinguished triangle in $\Db{B}$:
$$\xymatrix{\bar{Q}_V^+\ar[r]^{i_V} & F(V)\ar[r]^{\pi_V} &
\bar{F}(V)\ar[r]^{\alpha_V} & \bar{Q}_V^+[1].}$$

(1) For a morphism $f: V\lra V[k]$, we can form the following
commutative diagram in $\Db{B}$
$$\xymatrix{\bar{Q}_V^+\ar[r]^{i_V}\ar[d]^{a_f} &
F(V)\ar[r]^{\pi_V}\ar[d]^{F(f)} &
\bar{F}(V)\ar[r]^{\alpha_V}\ar[d]^{b_f} & \bar{Q}_V^+[1]\ar[d]^{a_f[1]}\\
\bar{Q}_V^+[k]\ar[r]^{i_V[k]} & F(V)[k]\ar[r]^{\pi_V[k]} &
\bar{F}(V)[k]\ar[r]^{\alpha_V[k]} & \bar{Q}_V^+[k+1]. }$$ The map
$b_f$ exists because the composition $i_VF(f)(\pi_V[k])$ belongs to
$\Hom_{\Db{B}}(\bar{Q}_V^+, \bar{F}(V)[k]) = 0$. If there is another
map $b_f': \bar{F}(V)\lra\bar{F}(V)[k]$ such that
$\pi_Vb'_f=F(f)(\pi_V[k])$, then $\pi_V(b_f-b'_f)=0$, and $b_f-b'_f$
factorizes through $\bar{Q}_V^+$. But $\Hom_{\Db{B}}(\bar{Q}_V^+[1],
\bar{F}(V)[k])\simeq \Hom_{\Kb{B}}(\bar{Q}_V^+[1], \bar{F}(V)[k])$ =
$0$. Hence $b_f=b'_f$, that is, the map $b_f$ is uniquely determined
by the above commutative diagram. Thus, we can define a morphism
$$\theta_k: \Hom_{\Db{A}}(V, V[k])\lra \Hom_{\Db{B}}(\bar{F}(V), \bar{F}(V)[k])$$
by sending $f$ to $b_f$. We claim that this $\theta_k$ is an
isomorphism.

In fact, it is surjective: For each map $b: \bar{F}(V)\lra
\bar{F}(V)[k]$, the composition $\pi_Vb(\alpha_V[k])$ belongs to
$\Hom_{\Db{B}}(F(V), \bar{Q}_V^+[k+1])\simeq \Hom_{\Db{A}}(GF(V),
G(\bar{Q}_V^+)[k+1])$. By Lemma \ref{AlmostV-property} (5), the
complex $G(\bar{Q}_V^+)$ is isomorphic in $\Db{A}$ to a bounded
complex $\cpx{P_V}$ of projective-injective $A$-modules such that
$P_V^i=0$ for all $i>1$. Hence $\Hom_{\Db{A}}(GF(V),
G(\bar{Q}_V^+)[k+1])$ $\simeq \Hom_{\Db{A}}(V, \cpx{P_V}[k+1])=0$,
and the map $\pi_Vb(\alpha_V[k])$ is zero. It follows that there is
a morphism $u: F(V)\lra F(V)[k]$ such that $u(\pi_V[k])=\pi_Vb$.
Since $F$ is an equivalence, we have $u=F(f)$ for some $f: V\lra
V[k]$, and consequently $b=\theta_k(f)$. This shows that $\theta_k$
is a surjective map.

Now we show that $\theta_k$ is injective: Assume that
$\theta_k(f)=b_f=0$. Then the composition $F(f)(\pi_V[k])=0$, and
consequently $F(f)$ factorizes through $\bar{Q}_V^+[k]$. It follows
that $GF(f)$ factorizes through $G(\bar{Q}_V^+)[k]\simeq
\cpx{P_V}[k]$, or equivalently, the map $f: V\lra V[k]$ factorizes
through the bounded complex $\cpx{P_V}$ of projective-injective
$A$-modules, say $f=gh$ for some $g: V\lra \cpx{P_V}$ and $h:
\cpx{P_V}\lra V[k]$. Since $k>0$, and since both $g$ and $h$ can be
chosen to be chain maps, we see immediately that $f=gh=0$. This
shows that the map $\theta_k$ is injective, and therefore $\theta_k$
is an isomorphism.

(2) By the above discussion, we have
$$\begin{array}{rl}
\pi_V \theta_{k}(g)(\theta_l(g)[k]) & =
(F(g)(\pi_V[k]))(\theta_l(h)[k])\\
& =F(g) ((\pi_V\theta_l(h))[k])\\
& = F(g) ((F(h)(\pi_V[l]))[k])\\
& =(F(g)F(h[k])) (\pi_V[k+l])\\
& = F(g(h[k]))(\pi_V[k+l]).
\end{array}$$
By the definition of $\theta_{k+l}$, we have
$\theta_{k+l}(g(h[k]))=\theta_{k}(g)(\theta_l(g)[k])$. $\square$

\medskip
{\bf\parindent=0pt Remark: } Let $f$ be in $\Hom_{\Db{A}}(V, V)$,
and let $g$ be in $\Hom_{\Db{A}}(V, V[k])$ for some $k>0$. If $t_f:
\bar{F}(V)\lra \bar{F}(V)$ is a morphism such that
$\pi_Vt_f=F(f)\pi_V$, then, by a proof similar to Lemma
\ref{thetamaps}(3), we have
$$t_f\theta_k(g)=\theta_k(fg)\mbox{\qquad and \qquad}
\theta_k(g)(t_f[k])=\theta_k(g(f[k])).$$

For instance, by Lemma \ref{kdiso}, we can assume that the map
$F(f): \cpx{\bar{Q}_V}\lra \cpx{\bar{Q}_V}$ is induced by a chain
map $\cpx{p}$, that is, $F(f)=\cpx{p}$ in $\Db{B}$. Since the map
$\pi_V$ is the canonical map from $\cpx{\bar{Q}_V}$ to
$\bar{Q}_V^0$, we see that the map $p^0: \bar{F}(V)\lra \bar{F}(V)$
satisfies the condition $\pi_Vp^0=F(f)\pi_V$. Therefore, by the
above discussion, we have
$$p^0\theta_k(g)=\theta_k(fg)\mbox{ and }
\theta_k(g)(p^0[k])=\theta_k(g(f[k])).$$

\medskip
\begin{Prop}
$\End_{\Kb{\pmodcat{\Ex{B}{\Phi}{N}}}}(\EP{B}{\Phi}{N}{\cpx{\bar{T}}})$
is isomorphic to $\Ex{A}{\Phi}{M}$. \label{PropEndYoneda}
\end{Prop}

\medskip
{\it Proof.} Let $(f_i)$ be in $\Ex{A}{\Phi}{M}$. By our assumption,
we have $\cpx{\bar{T}}=F(M)$. By Lemma \ref{kdiso}, the morphism
$F(f_0): \cpx{\bar{T}}\lra \cpx{\bar{T}}$ is equal in $\Db{B}$ to a
chain map. For simplicity, we shall assume that $F(f_0)$ is a chain
map. Recall that $\bar{F}(M)=\bar{T}^0$ by the definition of
$\bar{F}$ (see Lemma \ref{AlmostV-property} (3)\,).

Now we set $\Phi^+:=\Phi\backslash\{0\}$.   For each $k\in\Phi^+$,
by Lemma \ref{thetamaps}, we have a map $\theta_k(f_k):
\bar{F}(M)\lra \bar{F}(M)[k]$. This gives rise to a morphism
$$\mu\,\big(\iota_k(\theta_k(f_k))\big): \E{B}{\Phi}{N}{\bar{T}^0}\lra
\E{B}{\Phi}{N}{\bar{T}^0},$$ where $\mu$ is the isomorphism defined
in Lemma \ref{Extalgproperty} (1) and $\iota_k$ is the embedding
from $\Hom_{\Db{B}}(\bar{T}^0, \bar{T}^0[k])$ to
$\E{B}{\Phi}{\bar{T}^0}{\bar{T}^0}$. We claim that the composition
of $\mu\big(\iota_k(\theta_k(f_k))\big)$ with the differential
$\E{B}{\Phi}{N}{d}: \E{B}{\Phi}{N}{\bar{T}^0}$ $\ra
\E{B}{\Phi}{N}{\bar{T}^1}$ is zero.

Indeed, by the proof of Lemma \ref{Extalgproperty} (2), we have
$\E{B}{\Phi}{N}{d}=\mu\,(\iota_0(d))$. Thus,
$$\begin{array}{rll}
 \mu\,\big(\iota_k(\theta_k(f_k))\big)\E{B}{\Phi}{N}{d} &=
 \mu\,\big(\iota_k(\theta_k(f_k))\big)\mu\,(\iota_0(d)) & \quad\mbox{( by the proof of Lemma \ref{Extalgproperty} (2)\;)}\\
   & = \mu\,\big(\iota_k(\theta_k(f_k))\iota_0(d)\big) & \quad \mbox{( by Lemma \ref{Extalgproperty}(1)\;)}\\
   & =\mu\,\big(\iota_k(\theta_k(f_k)d[k])\big)=0 & \quad \mbox{( since }
   \theta_k(f_k)d[k]: \bar{T}^0\lra \bar{T}^1[k] \mbox{ must be zero\;).}
\end{array}$$

Thus, the map $\mu\big(\iota_k(\theta_k(f_k))\big)$ gives rise to an
endomorphism of $\EP{B}{\Phi}{N}{\cpx{\bar{T}}}$:
$$\xymatrix{
0\ar[r] &
\E{B}{\Phi}{N}{\bar{T}^0}\ar[r]^{\E{B}{\Phi}{N}{d}}\ar[d]_{\mu\,\big(\iota_k(\theta_k(f_k))\big)}
& \E{B}{\Phi}{N}{\bar{T}^1}\ar[r]\ar[d]^{0} &\cdots\ar[r]
&\E{B}{\Phi}{N}{\bar{T}^n}\ar[r]\ar[d]^{0} &
0\\
0\ar[r] & \E{B}{\Phi}{N}{\bar{T}^0}\ar[r]^{\E{B}{\Phi}{N}{d}} &
\E{B}{\Phi}{N}{\bar{T}^1}\ar[r] &\cdots\ar[r]
&\E{B}{\Phi}{N}{\bar{T}^n}\ar[r] & 0. }$$ We denote this
endomorphism by $\tilde{\theta}_k(f_k)$. Now, we define a map
$$\eta: \Ex{A}{\Phi}{M}\lra \End_{\Kb{\pmodcat{\Ex{B}{\Phi}{N}}}}(\EP{B}{\Phi}{N}{\cpx{\bar{T}}})$$
by sending $(f_i)$ to
$$\EP{B}{\Phi}{N}{F(f_0)}+\sum_{k\in\Phi^+}\tilde{\theta}_k(f_k).$$
We claim that $\eta$ is an algebra homomorphism. This will be shown
with help of the next lemma.

\begin{Lem}
  Let $(f_i)$ and $(g_i)$ be in $\Ex{A}{\Phi}{M}$, and let $k, l$
  be in $\Phi^+$. Then the following hold:

\medskip
$(1)$ $\tilde{\theta}_k(f_k)\tilde{\theta}_l(g_l)=\left\{
\begin{array}{ll}
\tilde{\theta}_{k+l}(f_k(g_l[k])), & {k+l\in\Phi ;} \\
0, & {k+l\not\in\Phi.}
\end{array}  \right.$

\medskip
$(2)$
$\EP{B}{\Phi}{N}{F(f_0)}\tilde{\theta}_k(g_k)=\tilde{\theta}_k(f_0g_k)$.

\medskip
$(3)$
$\tilde{\theta}_k(f_k)\EP{B}{\Phi}{N}{F(g_0)}=\tilde{\theta}_k(f_k(g_0[k]))$.\label{lemAlgHom}
\end{Lem}

\medskip
{\it Proof.} (1). By Lemma \ref{thetamaps} (2), we have
$$\iota_k\big(\theta_k(f_k)\big)\iota_{l}\big(\theta_l(f_l)\big)=\iota_{k+l}\big(\theta_{k+1}(f_k(g_l[k]))\big).$$
If $k+l\in\Phi$, then it follows that
$\tilde{\theta}_k(f_k)\tilde{\theta}_l(f_l)=\tilde{\theta}_{k+l}(f_k(g_l[k]))$
by applying $\mu$. If $k+l\not\in\Phi$, then $\iota_{k+l}=0$, and
consequently
$\iota_k\big(\theta_k(f_k)\big)\iota_{l}\big(\theta_l(f_l)\big)=0$.
Therefore $\tilde{\theta}_k(f_k)\tilde{\theta}_l(f_l)=0$ for
$k+l\not\in\Phi$.

(2) and (3). By definition, the map $\EP{B}{\Phi}{N}{F(f_0)}^0:
\E{B}{\Phi}{N}{\bar{T}^0}\lra \E{B}{\Phi}{N}{\bar{T}^0}$ is
$\E{B}{\Phi}{N}{F(f_0)^0}=\mu\,(\iota_0(F(f_0)^0))$, where
$F(f_0)^0: \bar{T}^0\lra\bar{T}^0$ is induced by the chain map
$F(f_0)$ from $\cpx{\bar{T}}$ to $\cpx{\bar{T}}$. By the remark just
before Lemma \ref{PropEndYoneda}, we have
$$\iota_0(F(f_0)^0)\iota_k\big(\theta_k(g_k)\big)=\iota_k\big(\theta_k(f_0g_k)\big) \mbox{ and  }
\iota_k\big(\theta_k(f_k)\big)\iota_0(F(g_0)^0)=\iota_k\big(\theta_k(f_k(g_0[k]))\big).$$
Applying $\mu$ to these equalities, one can easily see that
$$\EP{B}{\Phi}{N}{F(f_0)}\tilde{\theta}_k(g_k)=\tilde{\theta}_k(f_0g_k) \mbox{ and }
\tilde{\theta}_k(f_k)\EP{B}{\Phi}{N}{F(g_0)}=\tilde{\theta}_k(f_k(g_0[k])).$$
These are precisely the (2) and (3). $\square$

Now, we continue the proof of Lemma \ref{PropEndYoneda}: With Lemma
\ref{lemAlgHom} in hand, it is straightforward to check that $\eta$
is an algebra homomorphism. In the following we first show that
$\eta$ is injective.

Pick an $(f_i)$ in $\Ex{A}{\Phi}{M}$, let $\cpx{p}:=\eta((f_i))$.
Then we have
$$p^0=\E{B}{\Phi}{N}{F(f_0)^0}+\sum_{k\in\Phi^+}\mu\big(\iota_k(\theta_k(f_k))\big)$$
and $p^i=\E{B}{\Phi}{N}{F(f_0)^i}$ for all $i>0$. If $\cpx{p}=0$,
then there is map $h^i: \E{B}{\Phi}{N}{\bar{T}^i}\lra
\E{B}{\Phi}{N}{\bar{T}^{i-1}}$ for $i>0$ such that
$p^0=\E{B}{\Phi}{N}{d}h^1$ and
$p^i=\E{B}{\Phi}{N}{d}h^{i+1}+h^i\E{B}{\Phi}{N}{d}$ for all $i>0$.
Since $\bar{T}^i$ is projective-injective for all $i>0$, it follows
from Lemma \ref{Extalgproperty} (3) that, for each $i>0$, we have
$h^i=\E{B}{\Phi}{N}{u^i}$ for some $u^i:
\bar{T}^i\lra\bar{T}^{i-1}$. Hence
$$\E{B}{\Phi}{N}{F(f_0)^0}+\sum_{k\in\Phi^+}\mu\big(\iota_k(\theta_k(f_k))\big)=\E{B}{\Phi}{N}{d}\E{B}{\Phi}{N}{u^1}=\E{B}{\Phi}{N}{du^1}.$$
This yields that
$$\mu\big(\iota_0(F(f_0)^0-du^1)\big)=\E{B}{\Phi}{N}{F(f_0)^0-du^1}=\sum_{k\in\Phi^+}\mu\big(\iota_k(\theta_k(f_k))\big).$$
Since $\mu$ is an isomorphism, and since
$\E{B}{\Phi}{N}{\bar{T}^0}=\displaystyle{\bigoplus_{k\in\Phi}\,\Hom_{\Db{B}}(N,
\bar{T}^0[k])}$ is a direct sum, we get $F(f_0)^0=du^1$ and
$\theta_k(f_k)=0$ for all $k\in\Phi^+$. Since $\theta_k$ is an
isomorphism by Lemma \ref{thetamaps}, we have $f_k=0$ for all
$k\in\Phi^+$. Now for each $i>0$, we have
$$\E{B}{\Phi}{N}{F(f_0)^i}=p^i=\E{B}{\Phi}{N}{d}\E{B}{\Phi}{N}{u^{i+1}}+\E{B}{\Phi}{N}{u^i}\E{B}{\Phi}{N}{d}.$$
Hence $\E{B}{\Phi}{N}{F(f_0)^i-du^{i+1}-u^id}=0$. By Lemma
\ref{Extalgproperty} (2), the functor $\E{B}{\Phi}{N}{-}$ is
faithful. Therefore, we get $F(f_0)^i=du^{i+1}+u^id$ for $i>0$. Note
that we have shown that $F(f_0)^0=du^1$. Hence the morphism $F(f_0)$
is null-homotopic, that is, $F(f_0)=0$, and therefore $f_0=0$.
Altogether, we get $(f_i)=0$. This shows that $\eta$ is injective.

Finally, we show that $\eta$ is surjective. For $\cpx{p}$ in
$\End_{\Kb{\pmodcat{\Ex{B}{\Phi}{N}}}}(\EP{B}{\Phi}{N}{\cpx{\bar{T}}})$,
we can assume that $p^i=\E{B}{\Phi}{N}{t_i}$ with $t_i:
\bar{T}^i\lra \bar{T}^i$ for $i>0$ since $\bar{T}^i$ is
projective-injective for $i>0$. By Lemma \ref{Extalgproperty} (1),
we may assume further that
$p^0=\mu\,\big(\sum_{k\in\Phi}\iota_k(s_k)\big)$ with $s_k:
\bar{T}^0\lra\bar{T}^0[k]$ for $k\in\Phi$. By the proof of Lemma
\ref{Extalgproperty} (3), we have
$\mu\,\big(\iota_0(s_0)\big)=\E{B}{\Phi}{N}{s_0}$. Thus,
$p^0=\E{B}{\Phi}{N}{s_0}+\sum_{k\in\Phi^+}\mu\,\big(\iota_k(s_k)\big)$.
It follows from $\E{B}{\Phi}{N}{d}p^1=p^0\E{B}{\Phi}{N}{d}$ that
$$\displaystyle{\begin{array}{rl}
 \E{B}{\Phi}{N}{dt_1} & =\E{B}{\Phi}{N}{s_0d}+\displaystyle{\sum_{k\in\Phi^+}}\mu\,\big(\iota_k(s_k)\big)\mu\,\big(\iota_0(d)\big)\\
 & =\E{B}{\Phi}{N}{s_0d}+\displaystyle{\sum_{k\in\Phi^+}}\mu\,\big(\iota_k(s_k(d[k]))\big)\\
 & =\E{B}{\Phi}{N}{s_0d}   \quad (\mbox{because }s_k(d[k]): \bar{T}^0\lra \bar{T}^1[k]
 \mbox{ must be zero for } k>0).
 \end{array}}$$
Hence $dt_1=s_0d$ since $\E{B}{\Phi}{N}{-}$ is a faithful on
$\add(N)$. For each $i>0$, by the fact
$\E{B}{\Phi}{N}{d}p^{i+1}=p^i\E{B}{\Phi}{N}{d}$, we get
$dt_{i+1}=t_id$. This gives rise to a morphism $\cpx{\alpha}$ in
$\End_{\Kb{B}}(\cpx{\bar{T}})$ by defining $\alpha^0:=s_0$ and
$\alpha^i:=t_i$ for all $i>0$. By Lemma \ref{kdiso} and the fact
that $F$ is an equivalence, we conclude that $\cpx{\alpha}=F(f_0)$
for some $f_0\in\Hom_{\Db{A}}(M, M)$. The map
$\cpx{p}-\EP{B}{\Phi}{N}{\cpx{\alpha}}$ is a chain map $\cpx{\beta}$
from $\EP{B}{\Phi}{N}{\cpx{\bar{T}}}$ to itself with
$\beta^0=\sum_{k\in\Phi^+}\mu\,\big(\iota_k(s_k)\big)$ and
$\beta^k=0$ for all $k>0$. By Lemma \ref{thetamaps}, we can write
$s_k=\theta_k(f_k)$ with $f_k: M\lra M[k]$ for all $k\in\Phi^+$.
Thus
$\beta^0=\sum_{k\in\Phi^+}\mu\,\big(\iota_k(s_k)\big)=\sum_{k\in\Phi^+}\mu\,\big(\iota_k(\theta_k(f_k))\big)$,
and
$\cpx{p}-\EP{B}{\Phi}{N}{\cpx{\alpha}}=\sum_{k\in\Phi^+}\tilde{\theta}_k(f_k)$.
Consequently, we get
$$\cpx{p}=\EP{B}{\Phi}{N}{\cpx{\alpha}}+\sum_{k\in\Phi^+}\tilde{\theta}_k(f_k)
=\EP{B}{\Phi}{N}{F(f_0)}+\sum_{k\in\Phi^+}\tilde{\theta}_k(f_k)=\eta\big((f_i)\big)$$
for $(f_i)\in\Ex{A}{\Phi}{M}$. Hence $\eta$ is surjective. This
finishes the proof of Lemma \ref{PropEndYoneda}. $\square$

\begin{Lem}
Let $F: \Db{\Lambda}\lra\Db{\Gamma}$ be a derived equivalence
between Artin $R$-algebras $\Lambda$ and $\Gamma$, and let $\cpx{P}$
be a tilting complex associated to $F$. Suppose that the following
two conditions are satisfied.

\smallskip
$(1)$ All the terms of $\cpx{P}$ in negative degrees are zero, and
all the terms of $\cpx{P}$ in positive degrees are in
$\add(_{\Lambda}W)$ for some projective $\Lambda$-module
$_{\Lambda}W$ with $\add(\nu_{\Lambda}W)=\add({}_{\Lambda}W)$.

$(2)$ For the module $_{\Lambda}W$ in $(1)$, the complex
$F(_{\Lambda}W)$ is isomorphic to a complex in
$\Kb{\add(_{\Gamma}V)}$ for some projective $\Gamma$-module
$_{\Gamma}V$ with $\add(\nu_{\Gamma}V)=\add(_{\Gamma}V)$.

\smallskip
{\parindent=0pt Then} the quasi-inverse of $F$ is an almost
$\nu$-stable derived equivalence. \label{lemAlmostVstable}
\end{Lem}

\medskip
{\it Proof.} Let $G$ be a quasi-inverse of $F$. By the definition of
almost $\nu$-stable equivalences, we need to consider the tilting
complex associated to $G$. This is equivalent to considering
$F(\Lambda)$.

Since $\cpx{P}$ is a tilting complex over $\Lambda$, it is
well-known that $_{\Lambda}\Lambda$ is in
$\add(\bigoplus_{i\in\mathbb{Z}}P^i)$ which is contained in
$\add(P^0\oplus W)$ by the assumption (1). Hence
$F(_{\Lambda}\Lambda)$ is in $\add(F(P^0\oplus W))$. Let $P^+$ be
the complex $\sigma_{>0}\cpx{P}$. There is a distinguished triangle
$$\xymatrix{
P^+\ar[r] & \cpx{P}\ar[r] & P^0\ar[r] & P^+[1] }$$ in
$\Db{\Lambda}$. Applying $F$, we get a distinguished triangle
$$\xymatrix{
F(P^+)\ar[r] & F(\cpx{P})\ar[r] & F(P^0)\ar[r] & F(P^+)[1]
  }$$
in $\Db{\Gamma}$. By definition, there is an isomorphism
$F(\cpx{P})\simeq \Gamma$ in $\Db{\Gamma}$. By the assumption (1),
we have $P^+\in\Kb{\add(_{\Lambda}W)}$, and consequently $F(P^+)$ is
isomorphic in $\Db{\Gamma}$  to a complex $\cpx{R}$ in
$\Kb{\add(_{\Gamma}V)}$ by Assumption (2). Thus, the complex
$F(P^0)$ is isomorphic in $\Db{\Gamma}$ to the mapping cone of a
chain map from $\cpx{R}$ to $_{\Gamma}\Gamma$. This implies that
$F(P^0)$ is isomorphic in $\Db{\Gamma}$ to a complex $\cpx{S}$ in
$\Kb{\pmodcat{\Gamma}}$ such that $S^i\in\add(_{\Gamma}V)$ for all
$i\neq 0$. By the assumption (2) again, the complex $F(_{\Lambda}W)$
is isomorphic in $\Db{\Gamma}$ to a complex in
$\Kb{\add(_{\Gamma}V)}$. Hence $F(P^0)\oplus F(_{\Lambda}W)$ is
isomorphic in $\Db{\Gamma}$ to a complex $\cpx{U}$ in
$\Kb{\pmodcat{\Gamma}}$ such that $U^i\in\add(_{\Gamma}V)$ for all
$i\neq 0$. Note that $F(\Lambda)\in\add(F(P^0)\oplus
F(_{\Lambda}W))$. Therefore, the complex $F(\Lambda)$ is isomorphic
in $\Db{\Gamma}$ to a complex $\cpx{\bar{P}}$ in
$\Kb{\pmodcat{\Gamma}}$ such that $\bar{P}^i\in\add(_{\Gamma}V)$ for
all $i\neq 0$. Since $P^i=0$ for all $i<0$, we see from \cite[Lemma
2.1]{HX3} that $\cpx{\bar{P}}$ has zero homology in all positive
degrees. Hence we can assume that $\bar{P}^i=0$ for all $i>0$.

Thus, the complex $\cpx{\bar{P}}\simeq F(\Lambda)$ is a tilting
complex associated to $G$ and satisfies that $\bar{P}^i=0$ for all
$i>0$ and $\bar{P}^i\in\add(_{\Gamma}V)$ for all $i<0$. The complex
$\cpx{P}$ is a tilting complex associated to $F$ and satisfies that
$P^i=0$ for all $i<0$ and $P^i\in\add(_{\Lambda}W)$ for all $i>0$.
Since $\add(\nu_{\Lambda}W)=\add(_{\Lambda}W)$, and since
$\add(\nu_{\Gamma}V)=\add(_{\Gamma}V)$, it follows from
\cite[Proposition 3.8 (3)]{HX3} that the functor $G$ is an almost
$\nu$-stable derived equivalence. $\square$

Now we prove our main result, Theorem \ref{ThmYonedaAlgebra}, in
this section.

{\rm\bf Proof of Theorem \ref{ThmYonedaAlgebra}}. The statement (1)
follows from Lemma \ref{lemtiltYoneda}, Proposition
\ref{PropEndYoneda} and Lemma \ref{rickard}. It remains to prove
statement (2). Now we suppose that $\Phi$ is finite. Then
$\Ex{A}{\Phi}{M}$ and $\Ex{B}{\Phi}{N}$ are Artin $R$-algebras.

Let $_AE$ be a maximal $\nu$-stable $A$-module, and let $_B\bar{E}$
be a maximal $\nu$-stable $B$-module. Then $_AE$ can be viewed as a
direct summand of $_AM$. Let $\cpx{\bar{Q}_E}$ be $F(_AE)$ defined
in Lemma \ref{AlmostV-property} (1). Then $\cpx{\bar{Q}_E}$ is a
direct summand of
$\cpx{\bar{Q}_M}=\cpx{\bar{Q}}\oplus\cpx{\bar{Q}_X}$. Note that
$\cpx{\bar{Q}_M}$ is just the complex $\cpx{\bar{T}}$ considered in
Proposition \ref{PropEndYoneda}. Now we consider the isomorphism
$\eta$ in the proof of Proposition \ref{PropEndYoneda}. Let $e$ be
the idempotent in $\End_A(M)$ corresponding to the direct summand
$_AE$. Then $\iota_0(e)$ is the idempotent in $\Ex{A}{\Phi}{M}$
corresponding to the direct summand $\E{A}{\Phi}{M}{E}$ of
$\Ex{A}{\Phi}{M}$. By definition, $\eta\big(\iota_0(e)\big)$ is
$\EP{B}{\Phi}{N}{F(e)}$, which is the idempotent in
$\End_{\pmodcat{\Ex{B}{\Phi}{N}}}(\cpx{\bar{T}})$ corresponding to
$\EP{B}{\Phi}{N}{\cpx{\bar{Q}_E}}$. Hence the derived equivalence
$\hat{F}: \Db{\Ex{A}{\Phi}{M}}\lra \Db{\Ex{B}{\Phi}{N}}$ induces by
the isomorphism $\eta$ in the proof of Proposition
\ref{PropEndYoneda} sends $\E{A}{\Phi}{M}{E}$ to
$\EP{B}{\Phi}{N}{\cpx{\bar{Q}_E}}$. By \cite[Lemma 3.9]{HX3}, the
functor $F$ induces an equivalence between the triangulated
categories $\Kb{\add(_AE)}$ and $\Kb{\add(_B\bar{E})}$. Hence
$\cpx{\bar{Q}_E}=F(_AE)$ is in $\Kb{\add(_B\bar{E})}$, the complex
$\EP{B}{\Phi}{N}{\cpx{\bar{Q}_E}}$ belongs to
$\Kb{\add(\E{B}{\Phi}{N}{\bar{E}}}$ and consequently $\hat{F}$
induces a full, faithful triangle functor
$$\hat{F}: \Kb{\add(\E{A}{\Phi}{M}{E})}\lra \Kb{\add(\E{B}{\Phi}{N}{\bar{E}})}.$$
Since $\add(_AE)$ clearly generates $\Kb{\add(_AE)}$ as a
triangulated category, we see immediately that
$\add(\cpx{\bar{Q}_E})$ generates $\Kb{\add(_B\bar{E})}$ as a
triangulated category. This implies that
$\add(\EP{B}{\Phi}{N}{\cpx{\bar{Q}_E}})$ generates
$\Kb{\add(\E{B}{\Phi}{N}{\bar{E}})}$ as a triangulated category.
This shows that $$\hat{F}: \Kb{\add(\E{A}{\Phi}{M}{E})}\lra
\Kb{\add(\E{B}{\Phi}{N}{\bar{E}})}$$ is dense, and  therefore an
equivalence. Let $\hat{G}$ be a quasi-inverse of the derived
equivalence $\hat{F}$. Then the functor $\hat{G}$ also induces an
equivalence between the triangulated categories
$\Kb{\add(\E{B}{\Phi}{N}{\bar{E}})}$
 and $\Kb{\add(\E{A}{\Phi}{M}{E})}$. This implies that the complex $\hat{G}(\E{B}{\Phi}{N}{\bar{E}})$ is
isomorphic to a complex in $\Kb{\add(\E{A}{\Phi}{M}{E})}$.

Now we use Lemma \ref{lemAlmostVstable} to complete the proof. In
fact, the complex
 $\EP{B}{\Phi}{N}{\cpx{\bar{T}}}$ is a tilting complex associated to the
 derived equivalence $\hat{G}: \Db{\Ex{B}{\Phi}{N}}\lra \Db{\Ex{A}{\Phi}{M}}$. By
 definition, the $B$-module $\bar{Q}$ is in $\add(_B\bar{E})$.
 Thus, the term $\E{B}{\Phi}{N}{\bar{T}^i}$ of $\EP{B}{\Phi}{N}{\cpx{\bar{T}}}$
 in degree $i$ is in $\add(\E{B}{\Phi}{N}{\bar{E}})$ for all $i>0$, and it
 follows from Lemma \ref{Extalgproperty} (4) that
 $$\add(\nu_{\Ex{B}{\Phi}{N}}\E{B}{\Phi}{N}{\bar{E}})=
 \add(_{\Ex{B}{\Phi}{N}}\E{B}{\Phi}{N}{\nu_B\bar{E}})=\add(_{\Ex{B}{\Phi}{N}}\E{B}{\Phi}{N}{\bar{E}}).$$
Similarly, we have
$\add(\nu_{\Ex{A}{\Phi}{M}}\E{A}{\Phi}{M}{E})=\add(_{\Ex{A}{\Phi}{M}}\E{A}{\Phi}{M}{E})$.
Hence, by Lemma \ref{lemAlmostVstable}, the functor $\hat{F}$ is an
almost $\nu$-stable derived equivalence.

The statements on stable equivalence in Theorem
\ref{ThmYonedaAlgebra} follow from \cite[Theorem 1.1]{HX3}. This
finishes the proof. $\square$

\medskip
Note that the proof of Theorem \ref{ThmYonedaAlgebra} (2) shows also
that if both $\Ex{A}{\Phi}{M}$ and $\Ex{B}{\Phi}{N}$ are Artin
$R$-algebras, then the conclusion of Theorem \ref{ThmYonedaAlgebra}
(2) is valid.

Let us remark that, in case of finite-dimensional algebras over a
field, the special case for $\Phi=\Phi(1,0)$ in Theorem
\ref{ThmYonedaAlgebra} about stable equivalence was proved in
\cite[Proposition 6.1]{HX3} by using two-sided tilting complexes,
and the conclusion there guarantees a stable equivalence of Morita
type. But the proof there in \cite{HX3} does not work here any more,
since we do not have two-sided tilting complexes in general for
Artin algebras.

As a consequence of Theorem \ref{ThmYonedaAlgebra}, we have the
following corollary.

\begin{Koro}
Let $F:\Db{A}\lra\Db{B}$ be a derived equivalence between
self-injective Artin algebras $A$ and $B$, and let $\phi$ be the
stable equivalence induced by $F$. Then, for each $A$-module $X$ and
each admissible subset $\Phi$ of $\mathbb N$, the
$\Phi$-Auslander-Yoneda algebras $\Ex{A}{\Phi}{A\oplus X}$ and
$\Ex{B}{\Phi}{B\oplus \phi(X)}$ are derived-equivalent.
Particularly, the generalized Yoneda algebras $\Ext_A^*(A\oplus X)$
and $\Ext^*_B(B\oplus \phi(X))$ are derived-equivalent. Moreover, if
$\Phi$ is finite, then $\Ex{A}{\Phi}{A\oplus X}$ and
$\Ex{B}{\Phi}{B\oplus \phi(X)}$ are stably equivalent.
\label{twoselfinj}
\end{Koro}

{\it Proof.} There is an integer $i$ such that $F[i]$ is an almost
$\nu$-stable derived equivalence. Let $\phi_1$ be the stable
equivalence induced by $F[i]$. Then
$\phi(X)\simeq\phi_1\Omega^{i}(X)$ in $\stmodcat{B}$ for every
$A$-module $X$, where $\Omega^i$ is the $i$-th syzygy operator of
$A$. By the definition of an almost $\nu$-stable derived
equivalence, either $[i]$ or $[-i]$ is almost $\nu$-stable. Hence
$\Ex{}{\Phi}{A\oplus X}$ and $\Ex{A}{\Phi}{A\oplus\Omega^{i}(X)}$
are derived-equivalent by Theorem \ref{ThmYonedaAlgebra}. Thus, by
Theorem \ref{ThmYonedaAlgebra} again, the algebras
$\Ex{A}{\Phi}{A\oplus\Omega^{i}(X)}$ and
$\Ex{B}{\Phi}{B\oplus\phi_1\Omega^{i}(X)}$ are derived-equivalent.
The stable equivalence follows from \cite[Theorem 1.1]{HX3}. Thus
the proof is completed. $\square$

As a direct consequence of Corollary \ref{twoselfinj}, we have the
following corollary concerning Auslander algebras.

\begin{Koro}
Suppose that $A$ and $B$ are self-injective Artin algebras of finite
representation type. If $A$ and $B$ are derived-equivalent, then the
Auslander algebras of $A$ and $B$ are both derived and stably
equivalent.\label{DerivedAus}
\end{Koro}

Let us remark that the notion of a stable equivalence of Morita type
for finite-dimensional algebras can be formulated for Artin
$R$-algebras. But, in this case, we do not know if a stable
equivalence of Morita type between Artin algebras induces a stable
equivalence since we do not know whether a projective
$A$-$A$-bimodule is projective as a one-sided module when the ground
ring is a commutative Artin ring. So, Theorem \ref{ThmYonedaAlgebra}
(2), Corollary \ref{twoselfinj} (1) and Corollary \ref{DerivedAus}
ensure a stable equivalence between the endomorphism algebras of
generators over Artin algebras, while the main result in
\cite[Section 6]{HX3} ensures a stable equivalence of Morita type
between the endomorphism algebras of generators over
finite-dimensional algebras.

Note that if $A$ and $B$ are not self-injective, then Corollary
\ref{DerivedAus} may fail. For a counterexample, we just check the
following two algebras $A$ and $B$, where $A$ is given by the path
algebra of the quiver $\circ\ra \circ\ra \circ$, and $B$ is given by
$\circ\stackrel{\alpha}{\lra}\circ\stackrel{\beta}{\lra} \circ$ with
the relation $\alpha\beta=0$. Clearly, $B$ is the endomorphism
algebra of a tilting $A$-module. Note that the Auslander algebras of
$A$ and $B$ have different numbers of non-isomorphic simple modules,
and therefore are never derived-equivalent since derived
equivalences preserve the number of non-isomorphic simple modules
\cite{keller}. Notice that, though $A$ and $B$ are
derived-equivalent, there is no almost $\nu$-stable derived
equivalence between $A$ and $B$ since $A$ and $B$ are not stably
equivalent. This example shows also that Theorem
\ref{ThmYonedaAlgebra} may fail if we drop the almost $\nu$-stable
condition.

The following question relevant to Corollary \ref{DerivedAus} might
be of interest.

\medskip
{\parindent=0pt \bf Question.} Let $A$ and $B$ be self-injective
Artin algebras of finite representation type with $_AX$ and $_BY$
additive generators for $A$-mod and $B$-mod, respectively. Suppose
that there is a natural number $i$ such that the algebras
$\Ex{A}{\Phi(1,i)}{X}$ and $\Ex{B}{\Phi(1,i)}{Y}$ are
derived-equivalent. Are $A$ and $B$ derived-equivalent ?
\label{converseDeraus}

\medskip
We remark that Asashiba in \cite{asashiba} gave a complete
classification of representation-finite self-injective algebras up
to derived equivalence.

For a self-injective Artin $R$-algebra $A$, we know that the shift
functor [$-1$]: $\Db{A}\lra \Db{A}$ is an almost $\nu$-stable
derived equivalence, and this functor induces a stable functor
$\bar{F}: A$-\underline{mod}$\lra A$-\underline{mod}, which is
isomorphic to $\Omega_A(-)$, the Heller loop operator. Thus we have
the following corollary to Theorem \ref{ThmYonedaAlgebra}, which
extends \cite[Corollary 3.7]{hx2} in some sense.

\begin{Koro}
Let $A$ be a self-injective Artin algebra. Then, for any admissible
subset $\Phi$ of $\mathbb N$ and for any $A$-module $X$, we have a
derived equivalence between $\Ex{A}{\Phi}{A\oplus X}$ and
$\Ex{A}{\Phi}{A\oplus \Omega_A(X)}$. Moreover, if $\Phi$ is finite,
then there is an almost $\nu$-stable derived equivalence between
$\Ex{A}{\Phi}{A\oplus X}$ and $\Ex{A}{\Phi}{A\oplus \Omega_A(X)}$.
Thus they are stably equivalent.\label{selfinj}
\end{Koro}

Let us mention the following consequence of Corollary \ref{selfinj}.

\begin{Koro}
Let $A$ be a self-injective Artin algebra, and let $J$ be the
Jacobson radical of $A$ with the nilpotency index $n$. Then:

$(1)$ For any $1\le j\le n-1$ and for any admissible subset $\Phi$
of $\mathbb N$, the $\Phi$-Auslander-Yoneda algebras
$\Ex{A}{\Phi}{A\oplus \displaystyle{\bigoplus_{i=1}^{j}}A/J^i}$ and
$\Ex{A}{\Phi}{A\oplus \displaystyle{\bigoplus_{i=1}^{j}} J^i}$ are
derived-equivalent.

\smallskip
$(2)$ The global dimension of $\End_A(_AA\oplus J\oplus
J^2\oplus\cdots\oplus J^{n-1})$ is at most $n$.

\smallskip
$(3)$ The global dimension of $\End_A(_AA\oplus
\displaystyle{\bigoplus_{i= 1}^{n-1}}A/\soc^i(_AA))$ is at most $n$.

\smallskip
$(4)$ The global dimension of $\End_A(_AA\oplus \soc(_AA)\oplus
\cdots \oplus \soc^{n-1}(_AA))$ is at most
$n$.\label{auslanderconstruction}
\end{Koro}

{\it Proof.} Since the syzygy of $\bigoplus_{i=1}^{j}A/J^i$ is
$\bigoplus_{i=1}^{j} J^i$ up to a projective summand, we have (1)
immediately from Corollary \ref{selfinj}. The statement (2) follows
from \cite[Corollary 4.3]{HX3} together with a result of Auslander,
which says that, for any Artin algebra $A$, the global dimension of
$\End_A(A\oplus \bigoplus_{i=1}^{n-1}A/J^i)$ is at most $n$.

Since $_AA$ is injective, we know that $\add(_AA)=\add(D(A_A))$. It
follows from $D(A_A/J^i_A)\simeq \soc^i(D(A_A))$ that
$$\End_{A^{op}}(A_A\oplus \bigoplus_{i=1}^{n-1}A/J^i_A)\simeq
\Big(\End_A\big(D(A_A\oplus
\bigoplus_{i=1}^{n-1}A/J^i_A)\big)\Big)^{op}\simeq
\Big(\End_A\big(D(A_A)\oplus \bigoplus_{i=
1}^{n-1}\soc^i(D(A_A))\big)\Big)^{op}.$$ The latter is Morita
equivalent to $\big(\End_A(_AA\oplus \bigoplus_{i\ge
1}^{n-1}\soc^i(_AA))\big)^{op}$. This shows (4). The statement (3)
follows from (4), Corollary \ref{selfinj} and \cite[Corollary
4.3]{HX3}. $\square$

Finally, we state a dual version  of Theorem \ref{ThmYonedaAlgebra},
which will produce a derived equivalence between the endomorphism
algebras of cogenerators. First, we point out the following facts.

\begin{Lem}
Let $F:\Db{A}\lra\Db{B}$ be an almost $\nu$-stable derived
equivalence with a quasi-inverse functor $G$. Suppose $D$ is the
usual duality. Then we have the following.

$(1)$ The functor $DGD: \Db{B\opp}\lra\Db{A\opp}$ is an almost
$\nu$-stable derived equivalence with a quasi-inverse functor $DFD$.

$(2)$ Let $\bar{F}:\stmodcat{A}\lra \stmodcat{B}$ and
$\overline{DFD}:\stmodcat{A\opp}\lra\stmodcat{B\opp}$ be the stable
equivalence defined in {\rm Lemma \ref{AlmostV-property} (3) and
(4)}, respectively. Then,  for each $A$-module $X$, there is an
isomorphism $D\bar{F}(X)\simeq \overline{DFD}(D(X))$ in
$\modcat{B\opp}$.\label{propOpAlg}
\end{Lem}

{\it Proof.} (1)  Suppose that $\cpx{Q}$ and $\cpx{\bar{Q}}$ are
tilting complexes associated to $F$ and $G$, respectively. We assume
that $\cpx{Q}$ and $\cpx{\bar{Q}}$ are radical complexes. There is
an isomorphism
$$DGD(\HomP_B(\cpx{\bar{Q}}, {}_BB))\simeq DG(\nu_B\cpx{\bar{Q}})
\simeq D\nu_BG(\cpx{\bar{Q}})\simeq D\nu_A(_AA)\simeq \Hom_A({}_AA,
{}_AA)\simeq A_A$$ Similarly, we have $DFD(\HomP_A(\cpx{Q},
{}_AA))\simeq B_B$.
  Consequently, the complexes $\cpx{P}:=\HomP_B(\cpx{\bar{Q}}, {}_BB)$ and $\cpx{\bar{P}}:=\HomP_A(\cpx{Q},
{}_AA)$ are tilting complexes associated to $DGD$ and $DFD$,
respectively. Since $_B\bar{Q}=\bigoplus_{i=1}^n\bar{Q}^i$, we have
$\Hom_B(\bar{Q}, {}_BB)=\bigoplus_{i=1}^nP^{-i}$.  Moreover,
$$\nu_{B\opp}(\Hom_B(\bar{Q},{}_BB))\simeq D(_B\bar{Q})\simeq \Hom_B(\nu_B^-(\bar{Q}), {}_BB)\in\add(\Hom_B(\bar{Q},{}_BB))$$
 since $\nu_B^-\bar{Q}$ is in $\add(_B\bar{Q})$. Hence
 $\add(\nu_{B\opp}(\Hom_B(\bar{Q},{}_BB)))=\add(\Hom_B(\bar{Q},{}_BB))$.
Similarly, we have $\Hom_A({Q}, {}_AA)=\bigoplus_{i=1}^n\bar{P}^{i}$
and $\add(\nu_{A\opp}(\Hom_A(Q,{}_AA)))=\add(\Hom_A({Q},{}_AA))$,
and consequently $DGD$ is an almost $\nu$-stable derived
equivalence. Clearly, the functors $DGD$ and $DFD$ are mutually
quasi-inverse functors. This proves (1).

(2) For each $A$-module $X$, we have $DFD(D(X))=DF(X)$. By Lemma
\ref{AlmostV-property} (2), the complex $DFD(D(X))$ is isomorphic to
a complex $\cpx{P_{D(X)}}$ of the form
$$0\lra P_{D(X)}^{-n}\lra \cdots\lra P_{D(X)}^0\lra 0$$
with $P_{D(X)}^i\in\add(\Hom_B(\bar{Q}, {}_BB))$ for all $i<0$ and
$\overline{DFD}(D(X))=P_{D(X)}^0$. Consequently, the complex $F(X)$
is isomorphic to $D(\cpx{P_{D(X)}})$ of the form
$$0\lra D(P_{D(X)}^{0})\lra \cdots\lra D(P_{D(X)}^{-n})\lra 0$$
with $D(P_{D(X)}^{0})$ being in degree zero and
$D(P_{D(X)}^i)\in\add(\nu_B\bar{Q})=\add(_B\bar{Q})$  for all $i>0$.
By Lemma \ref{AlmostV-property} (1) and (3), we have
$\bar{F}(X)\simeq D(P_{D(X)}^0)=D\overline{DFD}(D(X))$ in
$\modcat{B}$. This finishes the proof. $\square$

Clearly, for an Artin algebra $A$ and an $A$-module $V$, the algebra
$\Ex{\Lambda}{\Phi}{V}$ is isomorphic to the opposite algebra of
$\Ex{\Lambda\opp}{\Phi}{D(V)}$ for every admissible subset $\Phi$ of
$\mathbb{N}$.

\begin{Koro}
Let $F:\Db{A}\lra \Db{B}$ be an almost $\nu$-stable derived
equivalence between two Artin algebras $A$ and $B$, and let
$\bar{F}$ be the stable equivalence defined in {\rm Lemma
\ref{AlmostV-property}}. For each $A$-module $X$, set
$M=D(A_A)\oplus X$ and $N=D(B_B)\oplus \bar{F}(X)$.  Suppose that
$\Phi$ is an admissible subset of $\mathbb{N}$. Then

$(1)$ The $\Phi$-Auslander-Yoneda algebras $\Ex{A}{\Phi}{M}$ and
$\Ex{B}{\Phi}{N}$ are derived-equivalent.

$(2)$ If $\Phi$ is finite, then there is an almost $\nu$-stale
derived equivalence between  $\Ex{A}{\Phi}{M}$ and
$\Ex{B}{\Phi}{N}$. \label{dual}
\end{Koro}

{\it Proof.}
 We consider the $A\opp$-module $DM=A_A\oplus D(X)$ and the $B\opp$-module $DN=B_B\oplus
 D\bar{F}(X)$.  By Lemma \ref{propOpAlg}, we see that
 $D\bar{F}(X)\simeq \overline{DFD}(D(X))$. Let $G$ be a
 quasi-inverse of $F$. Then the functor $DGD$ is an almost
 $\nu$-stable derived equivalence by Lemma \ref{propOpAlg} (1),
 and $\overline{DFD}$ is a quasi-inverse of $\overline{DGD}$. Thus,
 by Theorem \ref{ThmYonedaAlgebra} and by Lemma
 \ref{propOpAlg} (1), the corollary follows.
$\square$

\section{Derived equivalences for quotient algebras \label{sect4}}

In the previous section, we have seen  that there are many derived
equivalences between quotient algebras of $\Phi$-Auslander-Yoneda
algebras that are derived-equivalent (see Theorem
\ref{ThmYonedaAlgebra} and Subsection \ref{subsect3.1}). This
phenomenon gives rise to a general question: How to construct a new
derived equivalence for quotient algebras from the given one between
two given algebras ? In this section, we shall consider this
question and provide methods to transfer a derived equivalence
between two given algebras to a derived equivalence between their
quotient algebras. In particular, we shall prove Theorem \ref{theo2}

\subsection{Derived equivalences for algebras modulo ideals}

Let us start with the following general setting.

Suppose that $A$ is an Artin $R$-algebra over a commutative Artin
ring $R$, and suppose that $I$ is an ideal in $A$. We denote by
$\overline{A}$ the quotient algebra $A/I$ of $A$ by the ideal $I$.
The category $\overline{A}$-mod can be regarded as a full
subcategory of $A$-mod. Also, there is a canonical functor from
$A$-mod to $\overline{A}$-mod which sends each $X \in A$-mod to
$\overline{X}:= X/IX$. This functor induces a functor $^{-}:
\C{A}\lra \C{\overline{A}}$, which is defined as follows: for a
complex $\cpx{X} = (X^i)_{i\in {\mathbb Z}}$ of $A$-modules, let
$IX^{\bullet}$ be the sub-complex of $\cpx{X}$ in which the $i$-th
term is the submodule $IX^i$ of $X^i;$ we define
$\cpx{\overline{X}}$ to be the quotient complex of $\cpx{X}$ modulo
$IX^{\bullet}$. The action of $^{-}$ on a chain map can be defined
canonically. Thus $^{-}$ is a well-defined functor. For each complex
$\cpx{X}$ of $A$-modules, we have the following canonical exact
sequence of complexes:

$$\begin{CD}
0\lra I\cpx{X}\stackrel{\cpx{i}}{\lra}
\cpx{X}\stackrel{\cpx{\pi}}{\lra} \cpx{\overline{X}} \lra 0.
\end{CD}$$
For a complex $\cpx{Y}$ of $\overline{A}$-modules, this sequence
induces another exact sequence of $R$-modules:
$$\xymatrix{
0\ar[r] &\Hom_{\C{A}}(\cpx{\overline{X}}, \cpx{Y})\ar[r]^{{\pi}^*} &
\Hom_{\C{A}}(\cpx{X}, \cpx{Y})\ar[r]^{{i}^*} &
\Hom_{\C{A}}(I\cpx{X},\cpx{Y}).
 }$$
Since $\cpx{Y}$ is a complex of $\overline{A}$-modules, the map
$i^*$ must be zero, and consequently $\pi^*$ is an isomorphism. Now
we show that $\pi^*$ actually induces an isomorphism between
$\Hom_{\K{A}}(\cpx{\overline{X}}, \cpx{Y})$ and
$\Hom_{\K{A}}(\cpx{X}, \cpx{Y})$.

\begin{Lem}
Suppose that $A$ is an Artin algebra and $I$ is an ideal in $A$. Let
$\overline{A}$ be the quotient algebra of $A$ modulo $I$. If
$\cpx{X}$ is a complex of $A$-modules and $\cpx{Y}$ is a complex of
$\overline{A}$-modules, then we have a natural isomorphism of
$R$-modules
$$\pi^*: \Hom_{\K{A}}(\cpx{\overline{X}}, \cpx{Y})\longrightarrow
\Hom_{\K{A}}(\cpx{X}, \cpx{Y}).$$ \label{4.1}
\end{Lem}

{\it Proof.} Note that we have already an isomorphism
$$\pi^*: \Hom_{\C{A}}(\cpx{\overline{X}}, \cpx{Y})\longrightarrow
\Hom_{\C{A}}(\cpx{X}, \cpx{Y}).$$ Clearly, $\pi^*$ sends
null-homotopic maps to null-homotopic maps. This means that $\pi^*$
induces an epimorphism
$$\pi^*: \Hom_{\K{A}}(\cpx{\overline{X}}, \cpx{Y})\lra \Hom_{\K{A}}(\cpx{X}, \cpx{Y}). $$
Now let $\cpx{f}: \cpx{\overline{X}}\rightarrow \cpx{Y}$ be a chain
map such that $\pi^*(\cpx{f})=\cpx{\pi}\cpx{f}$ is null-homotopic.
Then there is a homomorphism $h^i: X^i\rightarrow Y^{i-1}$ for each
integer $i$ such that $\pi^if^i=h^id_Y^{i-1}+d_X^ih^{i+1}$. Note
that $h^i$ factorizes through $\pi^i$, that is, $h^i=\pi^i g^i$ for
some $g^i: \overline{X}^i\ra Y^{i-1}$. Hence we have
$$\begin{array}{rl}
  \pi^if^i & =h^id_Y^{i-1}+d_X^ih^{i+1}\\ &\\
           & =\pi^ig^id_Y^{i-1}+d_X^i\pi^{i+1}g^{i+1}\\ &\\
           & =\pi^ig^id_Y^{i-1}+\pi^id_{\overline{X}}^ig^{i+1}\\ &\\
           & =\pi^i(g^id_Y^{i-1}+d_{\overline{X}}^ig^{i+1}).\\
\end{array}$$
It follows that $f^i=g^id_Y^{i-1}+d_{\overline{X}}^ig^{i+1}$ since
$\pi^i$ is surjective for each $i$. Therefore, the map $\cpx{f}$ is
null-homotopic. Thus $\pi^*$ is injective. $\square$

\medskip
For any complexes $\cpx{X}$ and $\cpx{X'}$ over $\modcat{A}$, we
have a natural map
$$\eta: \Hom_{\K{A}}(\cpx{X},\cpx{X'})\lra
\Hom_{\K{A}}(\cpx{\overline{X}},\cpx{\overline{X'}}),$$ which is the
composition of $\cpx{\pi}_*: \Hom_{\K{A}}(\cpx{X}, \cpx{X'})\lra
\Hom_{\K{A}}(\cpx{X}, \cpx{\overline{X'}})$ with the map
$(\pi^*)^{-1}: \Hom_{\K{A}}(\cpx{X}, \cpx{\overline{X'}})\lra
\Hom_{\K{A}}(\cpx{\overline{X}}, \cpx{\overline{X'}})$ defined in
Lemma \ref{4.1}. In particular, if $\cpx{X}=\cpx{X'}$, then we get
an algebra homomorphism
$$\eta: \End_{\K{A}}(\cpx{X})\lra \End_{\K{A}}(\cpx{\overline{X}}).$$
Now, let $\cpx{T}$ be a tilting complex over $A$, and let
$B=\End_{\K{A}}(\cpx{T})$. Further, suppose that $I$ is an ideal in
$A$. By the above discussion, there is an algebra homomorphism
$$\eta: \End_{\K{A}}(\cpx{T})\lra \End_{\K{A}}(\cpx{\overline{T}}).$$
Let $J_I$ be the kernel of $\eta$, which is an ideal of $B$. Since
$(\pi^*)^{-1}$ is an isomorphism, we see that $J_I$ is the kernel of
the map $\cpx{\pi}_*: \End_{\K{A}}(\cpx{T})\lra
\Hom_{\K{A}}(\cpx{T},\cpx{\overline{T}})$. In fact, $J_I$ is also
the set of all endomorphisms of $\cpx{T}$ which factorize through
the embedding $I\cpx{T}\lra \cpx{T}$. We denote quotient algebra
$B/J_I$ by $\overline{B}$.

In the following, we study when the complex $\cpx{\overline{T}}$ is
a tilting complex over the quotient algebra $\overline{A}$ and
induces a derived equivalence between $\overline{A}$ and
$\overline{B}$. The following result supplies an answer to this
question.

\begin{Theo}
Let $A$ be an Artin algebra, and let $\cpx{T}$ be a tilting complex
over $A$ with the endomorphism algebra $B = \End_{\Kb{A}}(\cpx{T})$.
Suppose that $I$ is an ideal in $A$, and $\overline{A}:=A/I$. Let
$\overline{B}$ be the quotient algebra of $B$ modulo $J_I$. Then
$\cpx{\overline{T}}$ is a tilting complex over $\overline{A}$ and
induces a derived equivalence between $\overline{A}$ and
$\overline{B}$ if and only if $\Hom_{\Kb{A}}(\cpx{T},
I\cpx{T}[i])=0$ for all $i\neq 0$ and
$\Hom_{\Kb{A}}(\cpx{\overline{T}}, \cpx{\overline{T}}[-1])=0$.
\label{derivedquotient}
\end{Theo}

{\it Proof.} First, we assume $\Hom_{\Kb{A}}(\cpx{T},I\cpx{T}[i]) =
0$ for all $i\neq 0$ and $\Hom_{\Kb{A}}(\cpx{\overline{T}},
\cpx{\overline{T}}[-1])$ = $0$. Applying the functor
$\Hom_{\Db{A}}(\cpx{T},-)$ to the distinguished triangle
$$I\cpx{T}\stackrel{\cpx{i}}{\lra} \cpx{T}\stackrel{\cpx{\pi}}{\lra}
\cpx{\overline{T}} \lra I\cpx{T}[1],$$ for each integer $i$, we get
an exact sequence
$$\Hom_{\Db{A}}(\cpx{T}, \cpx{T}[i])\lra
\Hom_{\Db{A}}(\cpx{T}, \cpx{\overline{T}}[i])\lra
\Hom_{\Db{A}}(\cpx{T}, I\cpx{T}[i+1]),$$ which is isomorphic to the
exact sequence
$$(*)\qquad  \Hom_{\Kb{A}}(\cpx{T}, \cpx{T}[i])\lra
\Hom_{\Kb{A}}(\cpx{T}, \cpx{\overline{T}}[i])\lra
\Hom_{\Kb{A}}(\cpx{T}, I\cpx{T}[i+1]).$$
 Since the first and third terms of ($*$) are zero for $i\neq 0,
 -1$, the middle term $\Hom_{\Kb{A}}(\cpx{T},\cpx{\overline{T}}[i])$
 must be zero for $i\neq 0, -1$. Thus, taking our assumption into
 account, we have
 $$\begin{array}{rl}
 \Hom_{\Kb{\pmodcat{\overline{A}}}}(\cpx{\overline{T}}, \cpx{\overline{T}}[i])
 & \simeq \Hom_{\Kb{A}}(\cpx{\overline{T}},
 \cpx{\overline{T}}[i])
\\
  & \simeq  \Hom_{\Kb{A}}(\cpx{T}, \cpx{\overline{T}}[i])\\
 & = 0
 \end{array}$$
for all $i\neq 0$. Thus $\cpx{\overline{T}}$ is self-orthogonal in
$\Db{\overline{A}}$.

Note that the functor
$$(A/I)\otimesL_A-: \Kb{\pmodcat{A}}\longrightarrow \Kb{\pmodcat{\overline{A}}}$$
sends $\cpx{T}$ to $\cpx{\overline{T}}$. Let $\cal C$ be the full
triangulated subcategory of $\Kb{\pmodcat{\overline{A}}}$ generated
by $\add (\cpx{\overline{T}})$, and let $\cal D$ be a full
triangulated subcategory of $\Kb{\pmodcat{A}}$ consisting of those
$\cpx{X}$ for which $(A/I)\otimesL_A\cpx{X}$ belongs to $\cal C$.
Then $\cal D$ contains $\add(\cpx{T})$. Therefore ${\cal D} =
\Kb{\pmodcat{A}}$, and consequently $\add (\overline{A})$ is
contained in $\cal C$. Thus ${\cal C} =
\Kb{\pmodcat{\overline{A}}}$, and $\cpx{\overline{T}}$ is a tilting
complex over the quotient algebra $\overline{A}$. Since
$\Hom_{\Kb{A}}(\cpx{T}, I\cpx{T}[1])=0$, by the exact sequence
$(*)$, we have a surjective map $\cpx{\pi}_*: \Hom_{\Kb{A}}(\cpx{T},
\cpx{T})\ra \Hom_{\Kb{A}}(\cpx{T}, \cpx{\overline{T}})$. Therefore,
the algebra homomorphism $\eta:
\End_{\Kb{A}}(\cpx{T})\ra\End_{\Kb{A}}(\cpx{\overline{T}})$ is an
epimorphism. Hence $$\overline{B}=
\End_{\Kb{A}}(\cpx{T})/\Ker(\eta)\simeq
\End_{\Kb{A}}(\cpx{\overline{T}})\simeq
\End_{\Kb{\overline{A}}}(\cpx{\overline{T}}).$$ Consequently, the
tilting complex $\cpx{\overline{T}}$ induces a derived equivalence
between $\overline{A}$ and
 $\overline{B}$.

Conversely, we assume that $\cpx{\overline{T}}$ is a tilting complex
over $\overline{A}$ and induces a derived equivalence between
$\overline{A}$ and $\overline{B}$. Then
$\Hom_{\K{\overline{A}}}(\cpx{\overline{T}},\cpx{\overline{T}}[i])=0$
for all $i\neq 0$. Note that, for each integer $i$, we have an exact
sequence
$$(**) \qquad \Hom_{\Kb{A}}(\cpx{T}, \cpx{\overline{T}}[i-1]) \lra
\Hom_{\Kb{A}}(\cpx{T}, I\cpx{T}[i])\lra \Hom_{\Kb{A}}(\cpx{T},
\cpx{T}[i]).$$ Since $\Hom_{\Kb{A}}(\cpx{T},
\cpx{\overline{T}}[i-1])\simeq
\Hom_{\K{\overline{A}}}(\cpx{\overline{T}},\cpx{\overline{T}}[i-1])$
and since $\cpx{T}$ is self-orthogonal, the first and third terms of
($**$) are zero for $i\neq 0, 1$. It follows that
$\Hom_{\Kb{A}}(\cpx{T}, I\cpx{T}[i])=0$ for all $i\neq 0,1$. We
claim that $\Hom_{\Kb{A}}(\cpx{T}, I\cpx{T}[1])=0$. Indeed, we
consider the following exact sequence
$$\xymatrix@R=1mm@C=5mm{
\Hom_{\Kb{A}}(\cpx{T}, I\cpx{T}) \ar[r]^{\cpx{i}_*} &
\Hom_{\Kb{A}}(\cpx{T}, \cpx{T}) \ar[r]^{\cpx{\pi}_*} &
\Hom_{\Kb{A}}(\cpx{T}, \cpx{\overline{T}})
\ar[r] &\\
\Hom_{\Kb{A}}(\cpx{T}, I\cpx{T}[1])\ar[r] & \Hom_{\Kb{A}}(\cpx{T},
\cpx{T}[1])=0. }$$ Since the kernel of $\cpx{\pi}_*$ is $J_I$, the
image of $\cpx{\pi}_*$ is isomorphic to $\overline{B}$ as
$R$-modules. But we already have
$\overline{B}\simeq\End_{\Kb{A}}(\cpx{\overline{T}})$, which is
isomorphic to $\Hom_{\Kb{A}}(\cpx{T}, \cpx{\overline{T}})$ as an
$R$-module. Hence the map $\cpx{\pi}_*$ is surjective, and
$\Hom_{\Kb{A}}(\cpx{T}, I\cpx{T}[1])=0$. Clearly,
$\Hom_{\Kb{A}}(\cpx{\overline{T}},\cpx{\overline{T}}[-1])=0$.
Altogether, we have shown that $\Hom_{\Kb{A}}(\cpx{T},
I\cpx{T}[i])=0$ for all $i\neq 0$ and
$\Hom_{\Kb{A}}(\cpx{\overline{T}}$, $\cpx{\overline{T}}[-1])=0$.
This completes the proof of Theorem \ref{derivedquotient}. $\square$

\subsection{Derived equivalences for self-injective algebras modulo socles}

In the following, we shall use Theorem \ref{derivedquotient} to
prove our second main result in this paper. Let us first prove the
following lemma.

\begin{Lem}
Let $A$ be a self-injective basic algebra, and let $P$ be a direct
 summand of $_AA$.

$(1)$ If $J$ is an ideal of $A$ such that $_AJ\simeq {}_A\soc(P)$,
then $J=\soc(P)$.

$(2)$ If $\cpx{T}$ is a radical tilting complex over $A$ such that
the endomorphism algebra of $\cpx{T}$ is self-injective and basic,
then $T^i\simeq \nu_AT^i$ for all integers $i$. \label{lemJT}
\end{Lem}

{\it Proof.} (1) Let $e$ be the sum of the idempotents corresponding
to the simple direct summands of $\soc(P)$. By assumption, we have
$J\subseteq\soc(A)$ and $eJ=J$. Hence $J=eJ\subseteq
e(\soc(A))=\soc(P)$, and consequently $J=\soc(P)$.

(2) Let $B$ be the endomorphism algebra of $\cpx{T}$. Then there is
a derived equivalence $F:\Db{A}\lra\Db{B}$ such that
$F(\cpx{T})\simeq B$. Since $B$ is a self-injective basic algebra,
and since $F$ commutes with the Nakayama functor $\nu$, we have
$F(\nu_A\cpx{T})\simeq\nu_B F(\cpx{T})\simeq \nu_BB\simeq B\simeq
F(\cpx{T})$. Consequently, we have $\cpx{T}\simeq\nu_A\cpx{T}$ in
$\Db{A}$. Since $A$ is self-injective, we see that $\nu_A\cpx{T}$ is
also a complex in $\Kb{\pmodcat{A}}$. Hence
$\nu_A\cpx{T}\simeq\cpx{T}$ in $\Kb{\pmodcat{A}}$, and consequently
$\cpx{T}\simeq\nu_A\cpx{T}$ in $\Cb{A}$ since both $\cpx{T}$ and
$\nu_A\cpx{T}$ are radical complexes. Thus, the statement (2)
follows. $\square$

\begin{Theo}   Suppose that $A$ and $B$ are basic self-injective Artin
algebras, and that  $F:\Db{A}\lra\Db{B}$ is a derived equivalence. Let $P$ be
a direct summand of $_AA$, and let $P'$ be a direct summand of $_BB$
such that $F(\soc(P))$ is isomorphic to $\soc(P')$. Then the
quotient algebras $A/\soc(P)$ and $B/\soc(P')$ are
derived-equivalent. \label{thm4.3}
\end{Theo}

{\it Proof.} Since $A$ and $B$ are basic self-injective algebras,
$\soc(P)$ and $\soc(P')$ are ideals in $A$ and $B$, respectively. In
the following, we shall verify that the conditions of Theorem
\ref{derivedquotient} are satisfied by the ideal soc$(P)$ in $A$ and
the tilting complex $\cpx{T}$ associated to $F$.

Since $F(\soc(P))$ is isomorphic to $\soc(P')$, we can assume that
$P={\bigoplus_{i=1}^{s} P_i}$ and $P'={\bigoplus_{i=1}^{s} P'_i}$,
where $P$ and $P'$ are indecomposable such that $F(\soc(P_i))$ is
isomorphic to $\soc(P'_i)$ for all $i = 1, \cdots, s$. Let $D_i$ be
the endomorphism ring of $\soc(P_i)$, which is a division ring.
Since $F(\soc(P_i))\simeq\soc(P'_i)$, we see that $D_i$ is
isomorphic to $\End_B(\soc(P'_i))$. Note that a radical map $f:
M_1\ra M_2$ between two projective modules $M_1$ and $M_2$ has image
contained in $\rad(M_2)$. Since all the differential maps of
$\cpx{T}$ are radical maps, the image of $d_T^k$ is contained in
$\rad (T^{k+1})$ for all integers $k$. It follows that
$$\begin{array}{rcl}
\Hom_A(T^n,\soc(P_i)) & \simeq &
\Hom_{\Kb{A}}(\cpx{T}[n],\soc(P_i))\\
& \simeq &\Hom_{\Db{A}}(\cpx{T}[n], \soc(P_i))\\
& \simeq & \Hom_{\Db{B}}(B[n], \soc(P'_i))\\  & = &  0
\end{array}$$
for all $n\neq 0$. Hence, for each integer $n\ne 0$, the module
$\nu_A^{-1}P_i$ is not a direct summand of $T^n$. Since
$T^n\simeq\nu_A T^n$ (Lemma \ref{lemJT}(2)), we infer that $P_i$ is
not a direct summand of $T^n$ for all $n\ne 0$. Recall that
$\Hom_{\Db{A}}(\cpx{T}, \soc(P_i))\simeq
\Hom_{\Db{B}}(B,\soc(P'_i))\simeq\soc(P'_i)$ as $D_i^{\rm
op}$-modules. Since $B$ is basic, we see that $\soc(P'_i)$ is one-
dimensional over $D_i^{\rm op}$. Hence $\Hom_{\Db{A}}(\cpx{T},
\soc(P_i))$ is one-dimensional over $D_i^{\rm  op}$. It follows that
$\nu_A^{-1}P_i$ is a direct summand of $T^0$ with multiplicity 1.
Since $\nu_A T^0\simeq T^0$, we see that $P_i$ is a direct summand
of $T^0$ with multiplicity 1. Note that $\soc(P_i)X = 0$ for any
$A$-module $X$ if $P_i$ is not a direct summand of $X$. Hence
$\soc(P_i)\cpx{T}$ is isomorphic to the stalk complex $\soc(P_i)P_i
= \soc(P_i)$. Therefore
$$\Hom_{\Kb{A}}(\cpx{T},\soc(P)\cpx{T}[n])=\Hom_{\Kb{A}}(\cpx{T},
\oplus_{i=1}^{s}\soc(P_i)[n])=0$$ for all $n\ne 0$.

Let $\cpx{\overline{T}}$ be the quotient complex
$\cpx{T}/(\soc(P)\cpx{T})$. There is a canonical triangle in
$\Db{A}$:
$$ \soc(P)\cpx{T}\stackrel{\lambda}{\lra} \cpx{T}\lra
\cpx{\overline{T}}\lra (\soc(P)\cpx{T})[1].$$ Applying
$\Hom_{\Db{A}}(\cpx{T}, -)$ to this triangle, we have an exact
sequence of $B$-modules:
$$ 0\lra
\Hom_{\Db{A}}(\cpx{T},\cpx{\overline{T}}[-1])\lra
\Hom_{\Db{A}}(\cpx{T},\soc(P)\cpx{T})\stackrel{\lambda_*}{\lra}
\Hom_{\Db{A}}(\cpx{T},\cpx{T}).$$ We claim that $\lambda_*$ is a
monomorphism. Since $\soc(P)\cpx{T}$ is isomorphic to
$\bigoplus_{i=1}^{s}\soc(P_i)\cpx{T}$, the map $\lambda$ can be
written as $(\lambda_1,\cdots,\lambda_s)^{tr}$, where $\lambda_i:
\soc(P_i)\cpx{T}\rightarrow\cpx{T}$ is the canonical map, and where
$tr$ stands for the transpose of a matrix. Now we consider the
following commutative diagram of $B$-modules:
$$\xymatrix{
\Hom_{\Db{A}}(\cpx{T},\soc(P_i)\cpx{T})\ar[r]^(.55){(\lambda_i)_*}\ar[d]^{\simeq}
& \Hom_{\Db{A}}(\cpx{T},\cpx{T})\ar[d]^{\simeq}\\
\Hom_B(B,F(\soc(P_i)\cpx{T}))\ar[r]^(.60){F(\lambda_i)_*} &
\Hom_B(B,B).}$$ Since $\lambda_i\neq 0$, we see that $F(\lambda_i)$
is nonzero. Moreover, $F(\soc(P_i)\cpx{T})\simeq
F(\soc(P_i))\simeq\soc(P'_i)$. This implies that $F(\soc(P_i))$ is a
simple $B$-module for all $i$. Hence $F(\lambda_i)_*$ must be
injective. To show that $\lambda_*$ is injective, it suffices to
show that $F(\lambda)_*$ is injective. This is equivalent to proving
that $(F(\lambda_1)_*,\cdots,F(\lambda_s)_*)^{tr}$ is injective. For
this, we use induction on $s$. If $s=1$, the foregoing discussion
shows that this is true. Now we assume $s>1$. Then the kernel $K$ of
$(F(\lambda_1)_*,\cdots,F(\lambda_s)_*)^{tr}$ is the pull-back of
$(F(\lambda_1)_*,\cdots, F(\lambda_{s-1})_*)^{tr}$ and
$F(\lambda_s)_*$ both of which are monomorphisms by induction
hypothesis. Thus $K$ is isomorphic to a submodule of both
$\Hom_{\Db{A}}(\cpx{T},\bigoplus_{i=1}^{s-1}\soc(P))$ and
$\Hom_{\Db{A}}(\cpx{T},\soc(P_s))$. However, the $B$-modules
$\Hom_{\Db{A}}(\cpx{T},\soc(P_i))\simeq\soc(P'_i), i=1, \cdots, s, $
are pairwise  non-isomorphic simple $B$-modules since $B$ is basic.
This implies that $K=0$. Hence $\lambda_*$ is injective, and
therefore $\Hom_{\Db{A}}(\cpx{T},\cpx{\overline{T}}[-1])=0$. Since
$$\Hom_{\Kb{A}}(\cpx{\overline{T}},\cpx{\overline{T}}[-1])\simeq
\Hom_{\Kb{A}}(\cpx{{T}},\cpx{\overline{T}}[-1])\simeq
\Hom_{\Db{A}}(\cpx{T},\cpx{\overline{T}}[-1]),$$  it follows that
$\Hom_{\Kb{A}}(\cpx{\overline{T}},\cpx{\overline{T}}[-1])=0$. Hence
the complex $\cpx{T}$ and the ideal $\soc(P)$ satisfy all conditions
in Theorem \ref{derivedquotient}. Thus $A/\soc(P)$ and $B/J$ are
derived-equivalent, where $J$ is the ideal of $B$ consisting of maps
$b$ factorizing through the canonical map $\soc(P)\cpx{T}\lra
\cpx{T}$. Moreover, $J$ is isomorphic to
$\Hom_{\Kb{A}}(\cpx{T},\soc(P))$ as $B$-modules, and the latter is
isomorphic to $\soc(P')$. By Lemma \ref{lemJT} (1), we have
$J=\soc(P')$, and the theorem is proved. $\square$

\medskip
We give a criterion to judge when a derived equivalence satisfies
the condition in Theorem \ref{thm4.3}.

\begin{Prop} Let $\cpx{T}=(T^i,d^i)$ be a tilting complex associated to a derived
equivalence $F$ between self-injective basic Artin algebras $A$ and
$B$, and let $P$ be an indecomposable projective $A$-module. Suppose
we have the following:

$(1)$ $P\not\in\add(\nu_AT^i)$ for all $i\neq 0$;

$(2)$ the multiplicity of $P$ as a direct summand of $\nu_AT^0$ is
one.

{\parindent=0pt Let} $\cpx{T_P}$ be the indecomposable direct
summand of $\cpx{T}$ such that $P$ is a direct summand of
$\nu_A(T_P^0)$, and let $\bar{P}$ be the projective $B$-module
$\nu_B(\Hom_{\Kb{\pmodcat{A}}}(\cpx{T},\cpx{T_P}))$. Then
$F(\soc(_AP))\simeq \soc(_B\bar{P})$.
\end{Prop}

{\it Proof.} We know that the Nakayama functor sends $P$ to the
injective envelope of top$(_AP)$. From (1) it follows that
$\Hom_A(T^i, \soc(_AP))=0$ for all $i\neq 0$. Consequently,
$\Hom_{\Db{A}}(\cpx{T}$, $\soc(_AP)[i])=0$ for all $i\neq 0$. This
means that $F(\soc(_AP))$ is isomorphic in $\Db{B}$ to a $B$-module
$X$ that is indecomposable. Now we have the following isomorphisms:
$$\Hom_B(B, X)\simeq \Hom_{\Db{A}}(\cpx{T}, \soc(_AP))\simeq
\Hom_{\Db{A}}(\cpx{T_P}, \soc(_AP))\simeq\Hom_B(\nu_B^-\bar{P},X).$$
Hence $\soc(_B\bar{P})$ is the only simple $B$-module which occurs
as a composition factor of $X$. If $X$ were not simple, then we
would get a nonzero homomorphism $X\ra \top(X)\ra \soc(X)\ra X$,
which is not an isomorphism. This is a contradiction since
$\End_B(X)\simeq\End_{\Db{B}}(F(\soc(_AP))$ $\simeq
\End_A(\soc(_AP))$ is a division ring.  Hence $X$ is simple and
isomorphic to $\soc(_B\bar{P})$. This finishes the proof. $\square$

\subsection{Derived equivalences for algebras modulo annihilators}
Now, we turn to another construction for derived-equivalent quotient
algebras by using idempotent elements, which can be regarded as
another consequence of Theorem \ref{derivedquotient}.

\begin{Lem}\label{idealnablae} Let $e$ be an idempotent of
an Artin algebra $A$. Then there is a unique left ideal $I$ of $A$,
which is maximal with respect to the property $eI = 0$. Moreover,
the $I$ is an ideal of $A$. If, in addition, $\add(Ae)=
\add(D(eA))$, then $Ie=0$.
\end{Lem}

{\it Proof.} Note that such a left ideal $I$  in $A$ exists, and any
left ideal $L$ in $A$ with $eL=0$ is contained in $I$. Clearly, $I$
is a left ideal in $A$. We have to show that $I$ is a right ideal in
$A$. Let $x\in A$ and $a\in I$. Since the right multiplying with $x$
is a homomorphism $\varphi$ from $_AA$ to $_AA$, we see that the
image $\varphi(I)$ of $I$ under $\varphi$ is a left ideal in $A$.
Since $eI=0$, we have $\varphi(I)\subseteq I$, and $ax\in I$.

Suppose add$(Ae)=\add(D(eA))$. It follows from $$0= eI= \Hom_A(Ae,
I)\simeq \Hom_A(I,D(eA))$$ that Hom$_A(I,Ae)=0$. Clearly, the map
$\psi: I\ra Ae$ giving by $x\mapsto xe$ is a homomorphism from $I$
to $Ae$. Thus $\psi=0$ and $Ie=0$.
 $\square$

\medskip
Let $A$ be an  Artin algebra and $e$ an idempotent of $A$ such that
$\add(Ae)=\add(D(eA))$. By a result in \cite{IdemTilting}, there is
a tilting complex $\cpx{T}$ associated to $e$, which is defined in
the following way:  suppose $\varphi$ is a minimal right $\add
(Ae)$-approximation of $A$. Then we form the following complex:
$$ \cpx{T_f}: \quad 0\lra Q_1\stackrel{\varphi}{\lra} A\lra 0$$
with $A$ in degree zero. Let $\cpx{T_e}:=(Ae)[1]$. The tilting
complex $\cpx{T}$ associated with $e$ is defined to be the direct
sum of $\cpx{T_e}$ and $\cpx{T_f}$. Let $\lambda_e:
\cpx{T_e}\rightarrow\cpx{T}$ be the canonical inclusion and $p_e:
\cpx{T}\rightarrow\cpx{T_e}$ the canonical projection. Then
$\tilde{e}:=p_e\lambda_e$ is an idempotent in
$B:=\End_{\K{A}}(\cpx{T})$, which corresponds to the summand
$\cpx{T_e}$ of $\cpx{T}$. Thus, there is a derived equivalence $F:
\Db{A}\longrightarrow\Db{B}$, which sends $\cpx{T_e}$ to
$B\tilde{e}$, and $\cpx{T_f}$ to $B(1-\tilde{e})$. Let $\nabla(e)$
and $\nabla(\tilde{e})$ be the ideal $I$ of $A$ and $B$ defined by
$e$ and $\tilde{e}$ in Lemma \ref{idealnablae}, respectively. With
these notations in mind, we have the following proposition.

\begin{Prop}\label{A/nablae}   Let
$A$ be an Artin algebra and $e$ an idempotent element in $A$ such
that $\add(D(eA))=\add (Ae)$. Suppose that
$\cpx{T}=\cpx{T_e}\oplus\cpx{T_f}$ is the   tilting complex defined
by the idempotent $e$ and $B=\End_{\K{A}}(\cpx{T})$. Let $\tilde{e}$
be the idempotent element in $\End_{\K{A}}(\cpx{T})$ corresponding
to $\cpx{T_e}$. Then $A/\nabla(e)$ is derived-equivalent to
$B/\nabla(\tilde{e})$.
\end{Prop}

{\it Proof.} Let $F: \Db{A}\lra \Db{B}$ be the derived equivalence
given by the tilting complex $\cpx{T}$. Then
$F(\cpx{T_e})=F(Ae)[1])\simeq B\tilde{e}$.

\medskip
The complex $\nabla(e)\cpx{T}$ is isomorphic to $\nabla(e)$ because,
by Lemma \ref{idealnablae}, we have $\nabla(e)Ae=0$ and
$\nabla(e)\cpx{T}=\nabla(e)$, which is a complex with the only
non-zero term $\nabla(e)$ in degree zero. It is easy to see that
$\Hom_{\K{A}}(\cpx{T}, \nabla(e)[i])=0$ for all $i\neq 0$. Let
$\cpx{\overline{T}}$ be the quotient complex
$\cpx{T}/\nabla(e)\cpx{T}$. Then $\cpx{\overline{T}}$ is of the
following form:
$$\xymatrix{ 0\ar[r] & Ae\oplus
Q_1\ar[r]^(.55){\left[{0\atop\overline{\varphi}}\right]} &
\overline{A}\ar[r] & 0 },$$ where $\overline{A}=A/\nabla(e)$, and
where $\overline{\varphi}$ is the composition of $\varphi$ with the
canonical surjection from $A$ to $\overline{A}$. Since
$\Hom_{\K{A}}(\cpx{T}_f, \cpx{T}_e[-1])=0$, we get
$\Hom_A(\Coker(\varphi), Ae)=0$. Moreover, since
$\Coker(\overline{\varphi})$ is a quotient module of
$\Coker(\varphi)$, we have $\Hom_A(\Coker(\overline{\varphi}),
Ae)=0$. Thus
$$\Hom_{\K{A}}(\cpx{\overline{T}}, \cpx{\overline{T}}[-1])=0.$$ By
Theorem \ref{derivedquotient}, $\cpx{\overline{T}}$ is a tilting
complex over $A/\nabla(e)$, and $A/\nabla(e)$ is derived-equivalent
to $B/J$, where $J=\{\cpx{\alpha}\in\End_{\K{A}}(\cpx{T})\mid
\cpx{\alpha}\cpx{\pi}=0\}$, and where the map $\cpx{\pi}$ is the
canonical map from $\cpx{T}$ to $\cpx{\overline{T}}$. Note that
$\nabla(e)\cpx{T_{e}}=0$. This allows us to rewrite $\cpx{\pi}$ as
$$\xymatrix@R=10mm@C=10mm{
\cpx{T_e}\oplus\cpx{T_f}\ar[r]^(.40){\left[{1\,\,\,\,\,\;\;\,
0\,\,\,\,\atop \,\, \, 0\,\, \,\,\,\,\,\,\,\,\, \cpx{\pi_f}}\right]}
& \cpx{T_e}\oplus\cpx{T_f}/(\nabla(e)\cpx{T_f}). }$$

For any $\cpx{\alpha}\in J$, we can write $\cpx{\alpha}$ as
$$\xymatrix@R=10mm@C=10mm{
\cpx{T_e}\oplus\cpx{T_{f}}\ar[r]^{\left[{\cpx{\alpha_{11}}\,\,\,\cpx{\alpha_{12}}\atop
\cpx{\alpha_{21}}\,\,\,\cpx{\alpha_{22}}}\right]} &
\cpx{T_e}\oplus\cpx{T_{f}} .  }$$ Since $\cpx{\alpha}\cpx{\pi}=0$,
we have $\cpx{\alpha_{11}}=0=\cpx{\alpha_{21}}$ and
$\cpx{\alpha_{12}}\cpx{\pi_f}=0=\cpx{\alpha_{22}}\cpx{\pi_f}$. Hence
$\cpx{\alpha_{12}}:\cpx{T_e}\rightarrow \cpx{T_f}$ factorizes
through $\nabla(e)\cpx{T_f}=\nabla(e)$. But
$\Hom_{\K{A}}(\cpx{T_e},\nabla(e))=0$. This implies that
$\cpx{\alpha_{12}}=0$. Consequently, $J$ consists of maps
$\cpx{\alpha}$ of the form   $$\left[{0\quad\; 0\atop\,\,\,\,
0\,\,\,\, \,\,\;\cpx{\alpha_{22}}}\right]$$ with
$\cpx{\alpha_{22}}\cpx{\pi_f}=0$. Therefore $\tilde{e}J=0$ and
$J\subseteq\nabla(\tilde{e})$. By the proof of Theorem
\ref{derivedquotient}, we know that the quotient $B$-module $B/J$ is
isomorphic to $\Hom_{\K{A}}(\cpx{T}, \cpx{\overline{T}})$. Note that
we have a distinguished triangle
$$  A/\nabla(e)\lra \cpx{\overline{T}}\lra (Ae\oplus Q_1)[1]\lra
(A/\nabla(e))[1]$$ in $\K{A}$. Applying the functor
$\Hom_{\K{A}}(\cpx{T}, -)$ to this triangle, we get an exact
sequence
$$ \Hom_{\K{A}}(\cpx{T},
A/\nabla(e))\lra \Hom_{\K{A}}(\cpx{T}, \cpx{\overline{T}})\lra
\Hom_{\K{A}}(\cpx{T}, (Ae\oplus Q_1)[1]).$$ By the maximality of
$\nabla(e)$, the quotient $A/\nabla(e)$ has no submodule $X$ with
$eX=0$. Since $\varphi$ is a right $\add(Ae)$-approximation of $A$,
we have $e(\Coker(\varphi))=0$. It follows that
$\Hom_A(\Coker(\varphi)$, $A/\nabla(e))=0$. Hence we have
$\Hom_{\K{A}}(\cpx{T}, A/\nabla(e))=0$. Consequently,
$\Hom_{\K{A}}(\cpx{T}, \cpx{\overline{T}})$ can be embedded in
$\Hom_{\K{A}}(\cpx{T}$, $(Ae\oplus Q_1)[1])$, which is in $\add
(B\tilde{e})=\add(D(\tilde{e}B))$. This means that $J$ is the
maximal submodule of $B$ with $\tilde{e}J=0$. Hence
$J=\nabla(\tilde{e})$, and this finishes the proof. $\square$

\medskip
We point out that there is another type of construction by passing
derived equivalences between two given algebras to that between
their quotient algebras, namely, forming endomorphism algebras
first, and then passing to stable endomorphism algebras. For details
of this construction, we refer the reader to \cite[Corollary 1.2,
Corollary 1.3]{hx2}.

Now, we end this paper by two simple examples to illustrate our
results.

{\parindent=0pt\bf Example 1.} Let $k$ be a field, and let $A$ be a
 $k$-algebra given by the quiver
$$
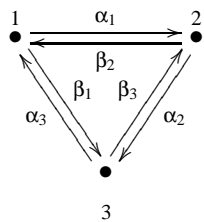
\xymatrix@R=6mm@C=8mm{
 \bullet\ar@<2pt>[rr]^{\alpha_1}^(0){1}^(1){2}\ar@<2pt>[rdd]^{\beta_1}
 & & \bullet\ar@<2pt>[ldd]^{\alpha_2}\ar@<2pt>[ll]^{\beta_2}\\ & & \\
 &\bullet\ar@<2pt>[luu]^(-.2){3}^{\alpha_3}\ar@<2pt>[ruu]^{\beta_3} &
} $$ with relations
$\alpha_i\beta_{i+1}-\beta_i\alpha_{i+2}=\alpha_i\alpha_{i+1}=\beta_i\beta_{i-1}=0$,
where the subscripts are considered modulo $3$. This algebra is
isomorphic to the group algebra of the alternative group $A_4$ if
$k$ has characteristic $2$. Let $e_2$ be the idempotent
corresponding to the vertex $2$, and let $\cpx{T}$ be the tilting
complex $\cpx{T}$ associated with $e_2$. Then the endomorphism
algebra $B$ of $\cpx{T}$ is given by the quiver
$$
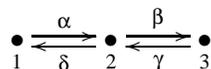
\xymatrix{
\bullet\ar@<2pt>[r]^{\alpha} &
\bullet\ar@<2pt>[r]^{\beta}\ar@<2pt>[l]^{\delta}^(1){1}
&\bullet\ar@<2pt>[l]^{\gamma}^(0){3}^(1){2}
 }$$
with relations $\alpha\delta=\gamma\beta=\delta\alpha\beta\gamma
-\beta\gamma\delta\alpha=0$. Note that $B$ is isomorphic to the
principal block of the group algebra of $A_5$ if $k$ has
characteristic $2$. It is easy to see that the idempotent
$\tilde{e}_2$ is the idempotent correspond to the vertex $2$ in the
quiver of $B$. Thus, by Proposition \ref{A/nablae}, the algebras
$A/\nabla(e_2)$ and $B/\nabla(\tilde{e}_2)$ are derived-equivalent.
A calculation shows that $A/\nabla(e_2)= A/\langle
\alpha_2\beta_3\rangle$ and $B/\nabla(\tilde{e}_2) =
B/\langle\beta\gamma\delta\alpha\rangle$. Note that the quotient
algebras $A/\langle\alpha_2\beta_3\rangle$ and
$B/\langle\beta\gamma\delta\alpha\rangle$ are stably equivalent of
Morita type by a result in \cite{LiuXi1}. Thus
$A/\langle\alpha_2\beta_3\rangle$ and
$B/\langle\beta\gamma\delta\alpha\rangle$ are not only
derived-equivalent, but also stably equivalent of Morita type.

\medskip
{\parindent=0pt\bf Example 2.} Let $m\ge 3$ be an integer, and let
$A=k[t]/(t^m)$, the quotient algebra of the polynomial algebra
$k[t]$ over a field $k$ in one variable $t$ modulo the ideal
generated by $t^m$. Let $X$ be the simple $A$-module $k$. Then
$\Ex{A}{\mathbb N}{A\oplus X}$ and $\Ex{A}{\mathbb N}{A\oplus
\Omega_A(X)}$ are infinite-dimensional $k$-algebras which can be
described by quivers with relations:
$$\begin{array}{ccc}
\Ex{A}{\mathbb N}{A\oplus X} & & \Ex{A}{\mathbb N}{A\oplus
\Omega_A(X)}\\ % & & \\
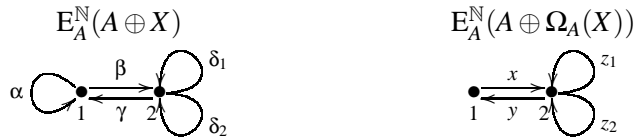
\xymatrix{ *{\bullet}\ar@<2pt>[r]^{\beta}\ar@(ul,dl)_{\alpha}
&*{\bullet}\ar@<2pt>[l]^(0.1){2}^(1){1}^{\gamma}\ar@(r,u)_(.4){\delta_1}\ar@(r,d)^(.4){\delta_2}
 } && \xymatrix{ *{\bullet}\ar@<2pt>[r]^{x}
&*{\bullet}\ar@<2pt>[l]^(0.1){2}^(1){1}^{y}\ar@(r,u)_(.4){z_1}\ar@(r,d)^(.4){z_2}
 }\\
&\qquad &\\
 \alpha^{m-1}-\beta\gamma=\alpha\beta=\gamma\alpha=\gamma\beta=0;
&& xz_i=z_iy=0, i=1,2;
\\
\delta_i\gamma=\beta\delta_i=0, i=1,2;&& z_1^2=z_1z_2-z_2z_1=0;\\
\delta_1^2=\delta_1\delta_2-\delta_2\delta_1=0. &&(yx)^{m-1}=0.
\end{array}$$
By Theorem \ref{ThmYonedaAlgebra}, or Corollary \ref{selfinj}, the
two algebras $\Ex{A}{\mathbb N}{A\oplus X}$ and $\Ex{A}{\mathbb
N}{A\oplus \Omega_A(X)}$ are derived-equivalent.

Let $n\ge 1$ be a natural number. Then the finite-dimensional
$k$-algebra $\Ex{A}{\Phi(1,n)}{A\oplus X}$ is the quotient of
$\Ex{A}{\mathbb N}{A\oplus X}$ by the ideal generated by
$\delta_2^{[\frac{n}{2}]+1}$ for $n$ an odd number, and by $
\delta_1\delta_2^{n/2}$ and $\delta_2^{n/2 +1}$ for $n$ an even
number, where $[\frac{n}{2}]$ is is the largest integer less than or
equal to $n/2$, and the finite-dimensional algebra
$\Ex{A}{\Phi(1,n)}{A\oplus \Omega_A(X)}$ is the quotient of
$\Ex{A}{\mathbb N}{A\oplus \Omega_A(X)}$ by the ideal generated by
$z_2^{[\frac{n}{2}]+1}$ for $n$ an odd number, and by $z_1z_2^{n/2}$
and $z_2^{n/2 +1}$ for $n$ an even number. By Corollary
\ref{selfinj}, we know that $\Ex{A}{\Phi(1,n)}{A\oplus X}$ and
$\Ex{A}{\Phi(1,n)}{A\oplus \Omega_A(X)}$  are both derived and
stably equivalent.

\medskip
{\bf Acknowledgements.} The corresponding author C.C.Xi thanks NSFC
(No.10731070) for limitedly partial support. The paper is revised
during Xi's visiting to RIMS at Kyoto University, he would like to
thank Prof. S.Ariki for his invitation and hospitality.

\medskip
April 28, 2009. Revised: November 11, 2009.

\end{document}